\documentclass[leqno,11pt]{amsart}
\usepackage{amsmath,amscd,amsthm,amsxtra}
\usepackage{epsfig,graphics,color,colortbl}
\usepackage{amssymb,latexsym}
\usepackage{mathrsfs,eucal,upgreek}
\usepackage[poly,all]{xy}
\usepackage{hyperref}
\usepackage{tikz}
\usetikzlibrary{decorations.pathreplacing}
\usetikzlibrary{snakes,arrows,shapes}
\usetikzlibrary{matrix,calc}
\usepackage{stackrel}
\usepackage{enumitem}
\newtheorem{thm}{\bf Theorem}[section]
\newtheorem{prop}[thm]{\bf Proposition}
\newtheorem{cor}[thm]{\bf Corollary}
\newtheorem{lem}[thm]{\bf Lemma}
\theoremstyle{definition}
\newtheorem{df}[thm]{Definition}
\newtheorem{rem}[thm]{Remark}
\newtheorem{exm}[thm]{\bf Example}

\numberwithin{equation}{section}

\newcommand{\A}{\mathcal{A}}

\newcommand{\B}{\mathbf{B}}

\newcommand{\W}{\mathcal{W}}
\newcommand{\cP}{\mathscr{P}}

\newcommand{\pf}{\noindent{\bfseries Proof. }}
\newcommand{\ov}{\overline}

\newcommand{\N}{\mathbb{N}}
\newcommand{\gl}{\mathfrak{gl}}
\newcommand{\Z}{\mathbb{Z}}

\newcommand{\te}{\widetilde{e}}
\newcommand{\tf}{\widetilde{f}}

\newcommand{\q}{\mathfrak{q}}

\newcommand{\mc}{\mathcal}
\newcommand{\mf}{\mathfrak}

\newcommand{\la}{\lambda}

\newcommand{\rev}{{\rm rev}}

\newcommand{\blue}[1]{{\color{blue}#1}}
\newcommand{\red}[1]{{\color{red}#1}}

\newcommand{\wt}{{\rm wt}}

\newcommand{\comp}{{\tt c}}
\newcommand{\compp}{{\tt c+}}
\newcommand{\rot}{{\tt \pi}}

\begin{document}
\title[Schur $P$-positive expansions]
{Crystals and Schur $P$-positive expansions}
\author{SEUNG-IL CHOI}
\address{Department of Mathematical Sciences, Seoul National University, Seoul 08826, Korea}
\email{ignatioschoi@snu.ac.kr}

\author{JAE-HOON KWON}
\address{Department of Mathematical Sciences, Seoul National University, Seoul 08826, Korea}
\email{jaehoonkw@snu.ac.kr}
\thanks{This work was supported by Samsung Science and Technology Foundation under Project Number SSTF-BA1501-01.}


\begin{abstract}
We give a new characterization of Littlewood-Richardson-Stembridge tableaux for Schur $P$-functions by using the theory of $\mf{q}(n)$-crystals. We also give alternate proofs of the Schur $P$-expansion of a skew Schur function due to Ardila and Serrano, and the Schur expansion of a Schur $P$-function due to Stembridge using the associated crystal structures. 
\end{abstract}
\maketitle

\section{Introduction}
Let $\cP^+$ be the set of strict partitions and let $P_\la$ be the Schur $P$-function corresponding to $\la\in \cP^+$ \cite{Sch}. The set of Schur $P$-functions is an important class of symmetric functions, which is closely related with representation theory and algebraic geometry (see \cite{Mac} and references therein). For example, the Schur $P$-polynomial $P_\la(x_1,\ldots,x_n)$ in $n$ variables is the character of a finite-dimensional irreducible representation $V_n(\la)$ of the queer Lie superalgebra $\mf{q}(n)$ with highest weight $\la$ up to a power of $2$ when the length $\ell(\la)$ of $\la$ is no more than $n$ \cite{Srgv85}. 

The set of Schur $P$-functions forms a basis of a subring of the ring of symmetric functions, and the structure constants with respect to this basis are non-negative integers, that is, given $\la, \mu, \nu \in \cP^+$,
\begin{equation*}
P_{\mu}  P_{\nu} 
=\underset{\la}{\sum}f^{\la}_{\mu \nu}P_{\la},
\end{equation*}
for some non-negative integers $f^{\la}_{\mu \nu}$. 
The first and the most well-known result on a combinatorial description of $f^{\la}_{\mu \nu}$ was obtained by Stembridge \cite{Ste} using shifted Young tableaux, which is a combinatorial model for Schur $P$- or $Q$-functions \cite{Sa,Wo}. 
It is shown that $f^\la_{\mu\nu}$ is equal to the number of semistandard tableaux with entries in a $\Z_2$-graded set ${\mc N}=\{\,1'<1<2'<2<\cdots \,\}$ of shifted skew shape $\la/\mu$ and weight $\nu$ such that (i) for each integer $k\geq 1$ the southwesternmost entry with value $k$ is unprimed or of even degree and (ii) the reading words satisfy the {\em lattice property}. 
Here we say that the value $|x|$ is $k$
when $x$ is either $k$ or $k'$ in a tableau.
Let us call these tableaux the {\em Littlewood-Richardson-Stembridge (LRS) tableaux} (Definitions \ref{df:St lattice property} and \ref{Def:Stembridge's description}).

Recently, two more descriptions of $f^\la_{\mu\nu}$ were obtained in terms of semistandard decomposition tableaux, which is another combinatorial model for Schur $P$-functions introduced by Serrano~\cite{Se}.  
It is shown by Cho that $f^\la_{\mu\nu}$ is given by the number of semistandard decomposition tableaux of shifted shape $\mu$ and weight $w_0(\la-\nu)$ whose reading words satisfy the {\em $\la$-good property} (see~\cite[Corollary 5.14]{Cho}).
Here we assume that $\ell(\la), \ell(\mu), \ell(\nu)\leq n$, and  $w_0$ denotes the longest element in the symmetric group $\mf S_n$.
Another description is given by  Grantcharov, Jung, Kang, Kashiwara, and Kim~\cite{GJKKK14}
based on their crystal base theory for the quantized enveloping algebra of $\mf{q}(n)$ \cite{GJKKK15}. 
They realize the crystal $\B_n(\la)$ associated to $V_n(\la)$ as the set of semistandard decomposition tableaux of shape $\la$ with entries in  $\{\,1 < 2 < \cdots < n\,\}$,
and describe $f^\la_{\mu\nu}$ by characterizing the lowest weight vectors of weight $w_0\la$ in the tensor product $\B_n(\mu)\otimes \B_n(\nu)$.
We also remark that bijections between the above mentioned combinatorial models for $f^\la_{\mu\nu}$ are studied in \cite{CNO14} using insertion schemes for semistandard decomposition tableaux.

The main result in this paper is to give another description of $f^\la_{\mu\nu}$ using the theory of $\mf{q}(n)$-crystals, and show that it is indeed equivalent to that of Stembridge. 
More precisely, we show that $f^\la_{\mu\nu}$ is equal to the number of semistandard tableaux with entries in $\mc N$ of shifted skew shape $\la/\mu$ and weight $\nu$
such that (i) for each integer $k\geq 1$ the southwesternmost entry with value $k$ is unprimed or of even degree  and 
(ii) the reading words satisfy the {\em ``lattice property"} 
(see Definitions \ref{Def:new lattice property} and \ref{Def:new LR tab} and Theorem \ref{Thm:our-characterization}). It is obtained by semistandardizing the standard tableaux which parametrize the lowest weight vectors counting $f^\la_{\mu\nu}$ in \cite{GJKKK14}, where the ``lattice property" naturally arises from the configuration of entries in semistandard decomposition tableaux. We show that these tableaux for $f^{\la}_{\mu\nu}$ are equal to LRS tableaux 
(Theorem \ref{thm:F=LRS}), and hence obtain a new characterization of LRS tableaux.

We study other Schur $P$- or $Q$-positive expansions and their combinatorial descriptions from a viewpoint of crystals.
First we consider the Schur $P$-positive expansion of a skew Schur function
\begin{equation*}
s_{\lambda/\delta_r} = \sum_{\nu\in\cP^+} a_{\lambda/\delta_r\,\nu} \hskip 0.5mm P_\nu
\end{equation*}
for a skew diagram $\lambda/\delta_r$ contained in a rectangle $((r+1)^{r+1})$, where $\delta_r=(r,r-1,\ldots,1)$ \cite{AS}. We give a combinatorial description of $a_{\lambda/\delta_r\,\nu}$ (Theorem \ref{Thm:Main-sect5}) by considering a $\mf{q}(n)$-crystal structure on the set of usual semistandard tableaux of shape $\la/\delta_r$ and characterizing the lowest weight vectors corresponding to each $\nu\in \cP^+$.
As a byproduct we also give a simple alternate proof of Ardila-Serrano's description of $a_{\lambda/\delta_r\,\nu}$ \cite{AS} (Theorem \ref{Thm:second-sect5}), which can be viewed as a standardization of our description. 

We next consider the Schur expansion of a Schur $P$-function
\begin{equation*}
P_\lambda = \sum_{\mu} g_{\lambda \mu} s_\mu
\end{equation*}
for $\la\in\cP^+$. It is equivalent to the expansion of a symmetric function $S_\mu=S_\mu(x,x)$ in terms of Schur $Q$-functions $Q_\la=2^{\ell(\la)}P_\la$, where $S_\mu(x,y)$ is a super Schur function in variables $x$ and $y$.
We give a simple and alternate proof of Stembridge's description of $g_{\la\mu}$ \cite{Ste} (Theorem \ref{thm-coeff-g}) by characterizing the type $A$ lowest weight vectors of weight $w_0\mu$ in the $\mf{q}(n)$-crystal $\B_n(\la)$ when $\ell(\la), \ell(\mu)\leq n$.

Finally, we introduce the notion of semistandard decomposition tableaux of shifted skew  shape. We consider a ${\mf q}(n)$-crystal structure on the set of such tableaux, and describe its decomposition into $\B_n(\la)$'s, which implies that the corresponding character has a Schur $P$-positive expansion though it is not equal to a skew Schur $P$-function in general.

The paper is organized as follows. In Section 2, we review the notion of $\mf{q}(n)$-crystals and related results. In Section 3, we describe a combinatorial description of $f^\la_{\mu\nu}$ and show that it is equivalent to that of Stembridge. In Sections 4 and 5, we discuss the Schur $P$-positive expansion of a skew Schur function and the Schur expansion of a Schur $P$-function, respectively. In Section 6, we discuss semistandard decomposition tableaux of shifted  skew shape, and the Schur $P$-positive expansions of their characters.

\section{Crystals for queer Lie superalgebras}
\subsection{Notation and terminology}  
In this subsection, we introduce necessary notations and terminologies.
Let $\Z_+$ be the set of non-negative integers.
We fix a positive integer $n\geq 2$ throughout this paper.

Let $\mathscr{P}=\{\,\la=(\la_i)_{i\geq 1}\,|\,\la_i\in \Z_+,\, \la_i\geq \la_{i+1}\, (i\geq 1),\, \sum_{i\geq 1}\la_i<\infty \,\}$ be the set of partitions, and let $\cP^+=\{\,\la=(\la_i)_{i\geq 1}\,|\,\la\in \cP,\, \la_i=\la_{i+1} \Rightarrow \la_i=0 \ (i\geq 1) \,\}$ be the set of strict partitions. 
For $\la\in \cP$, let $\ell(\la)$ denote the length of $\la$, and $|\la|=\sum_{i\geq 1}\la_i$.
Let $\cP_n=\{\,\la \,|\,\ell(\la)\leq n\,\}\subseteq \cP$ and $\cP^+_n=\cP^+\cap \cP_n$.

The (unshifted) diagram of $\lambda\in\cP$ is defined to be the set
\begin{displaymath}
D_\lambda = \{\, (i,j) \in \mathbb{N}^2 \,:\, 1 \leq j \leq \lambda_i, \ 1 \leq i \leq \ell(\lambda)\, \},
\end{displaymath}
and  
the shifted diagram of $\lambda\in\cP^+$ is defined to be the set
\begin{displaymath}
D^+_\lambda = \{\, (i,j)  \in \mathbb{N}^2 \,:\, i \leq j \leq \lambda_i+i-1, \ 1 \leq i \leq \ell(\lambda) \, \}.
\end{displaymath}
We identify $D_\la$ and $D^+_{\la}$ with diagrams where a box is placed at the $i$-th row from the top and the $j$-th column from the left for each $(i,j)\in D_{\la}$ and $D^+_{\la}$, respectively.
For instance, if $\lambda = (6,4,2,1)$, then 
\begin{displaymath}
\begin{tikzpicture}[baseline=0mm]
\def \hhh{3.5mm}    
\def \vvv{4.2mm}    
\node[left] at (-\hhh*0,-\vvv*0.5) {$D_{\lambda}=$};
\draw[-] (\hhh*0,0) rectangle (\hhh*1,\vvv*1);
\draw[-] (\hhh*1,0) rectangle (\hhh*2,\vvv*1);
\draw[-] (\hhh*2,0) rectangle (\hhh*3,\vvv*1);
\draw[-] (\hhh*3,0) rectangle (\hhh*4,\vvv*1);
\draw[-] (\hhh*4,0) rectangle (\hhh*5,\vvv*1);
\draw[-] (\hhh*5,0) rectangle (\hhh*6,\vvv*1);
\draw[-] (\hhh*0,-\vvv*1) rectangle (\hhh*1,\vvv*0);
\draw[-] (\hhh*1,-\vvv*1) rectangle (\hhh*2,\vvv*0);
\draw[-] (\hhh*2,-\vvv*1) rectangle (\hhh*3,\vvv*0);
\draw[-] (\hhh*3,-\vvv*1) rectangle (\hhh*4,\vvv*0);
\draw[-] (\hhh*0,-\vvv*2) rectangle (\hhh*1,-\vvv*1);
\draw[-] (\hhh*1,-\vvv*2) rectangle (\hhh*2,-\vvv*1);
\draw[-] (\hhh*0,-\vvv*3) rectangle (\hhh*1,-\vvv*2);
\end{tikzpicture}
\hskip 8mm
\begin{tikzpicture}[baseline=0mm]
\def \hhh{3.5mm}    
\def \vvv{4.2mm}    
\node[left] at (-\hhh*0,-\vvv*0.5) {and};
\end{tikzpicture}
\hskip 8mm
\begin{tikzpicture}[baseline=0mm]
\def \hhh{3.5mm}    
\def \vvv{4.2mm}    
\node[left] at (-\hhh*0,-\vvv*0.5) {$D^+_{\lambda}=$};
\draw[-] (\hhh*0,0) rectangle (\hhh*1,\vvv*1);
\draw[-] (\hhh*1,0) rectangle (\hhh*2,\vvv*1);
\draw[-] (\hhh*2,0) rectangle (\hhh*3,\vvv*1);
\draw[-] (\hhh*3,0) rectangle (\hhh*4,\vvv*1);
\draw[-] (\hhh*4,0) rectangle (\hhh*5,\vvv*1);
\draw[-] (\hhh*5,0) rectangle (\hhh*6,\vvv*1);
\draw[-] (\hhh*1,-\vvv*1) rectangle (\hhh*2,\vvv*0);
\draw[-] (\hhh*2,-\vvv*1) rectangle (\hhh*3,\vvv*0);
\draw[-] (\hhh*3,-\vvv*1) rectangle (\hhh*4,\vvv*0);
\draw[-] (\hhh*4,-\vvv*1) rectangle (\hhh*5,\vvv*0);
\draw[-] (\hhh*2,-\vvv*2) rectangle (\hhh*3,-\vvv*1);
\draw[-] (\hhh*3,-\vvv*2) rectangle (\hhh*4,-\vvv*1);
\draw[-] (\hhh*3,-\vvv*3) rectangle (\hhh*4,-\vvv*2);
\node at  (\hhh*6.3,-\vvv*3) {.};
\end{tikzpicture}
\end{displaymath} 

Let $\mc{A}$ be a linearly ordered set. We denote by $\W_\A$ the set of words of finite length with letters in $\A$. For $w\in \W_\A$ and $a\in \A$, let $c_a(w)$ be the number of occurrences of $a$ in $w$.

For $\la, \mu\in\cP$ with $D_\mu \subseteq D_\lambda$, 
a {\em tableau of shape $\lambda / \mu$}  means a filling on the skew diagram $D_\lambda \setminus D_\mu$ with entries in  $\mc{A}$.
For $\la, \mu \in \cP^+$ with $D^+_\mu \subseteq D^+_\la$, 
a {\em tableau of shifted  shape $\la/\mu$} is defined in a similar way. For a tableau $T$ of (shifted) shape $\la/\mu$, let $w(T)$ be the word given by reading the entries of $T$ row by row from top to bottom, and from right to left in each row. 
We denote by $T_{i,j}$ the $j$-th entry (from the left) of the $i$-th row of $T$ from the top. 
For $1\leq i\leq \ell(\la)$, let $T^{(i)}=T_{i,{\la_i}}\cdots T_{i,1}$ be the subword of $w(T)$ corresponding to the $i$-th row of $T$. 
Then we have $w(T)=T^{(1)}\cdots T^{(\ell(\la))}$. We denote by $w_{\rev}(T)$ the reverse word of $w(T)$. 
Note that $T_{i,j}$ is not the entry of $T$ at the $(i,j)$-position of the (shifted) skew diagram of $\la/\mu$, that is, $(i,j)\in D_{\la}\setminus D_{\mu}$ or $(i,j)\in D^+_{\la}\setminus D^+_{\mu}$. For $a\in\A$, let $c_a(T)=c_a(w(T))$ be the number of occurrences of $a$ in $T$.

Suppose that $\mc{A}$ is a linearly ordered  set with a $\mathbb{Z}_2$-grading $\mc{A}=\mc{A}_0\sqcup\mc{A}_1$. 
For $\la, \mu\in\cP$ with $D_\mu \subseteq D_\lambda$, let $SST_\A(\lambda/\mu)$ be the set of tableaux of shape $\lambda/\mu$ with entries in $\A$ which is semistandard, that is, (i) the entries in each row (resp. column) are
weakly increasing from left to right (resp. from top to bottom), 
(ii) the entries in $\mc{A}_0$ (resp. $\mc{A}_1$) are strictly increasing in each
column (resp. row). 
Similarly, for $\la, \mu\in \cP^+$ with $D^+_\mu \subseteq D^+_\la$, 
we define $SST^+_{\A}(\la/\mu)$ to be the set of semistandard tableaux of shifted shape $\lambda/\mu$  with entries in $\A$. 
 
Let ${\mc N}=\{\,1' < 1 < 2' < 2 <  \cdots \,\}$ be a linearly ordered set with a $\Z_2$-grading ${\mc N}_0=\N$ and ${\mc N}_1=\N'=\{1',2',\cdots\}$. Put $[n]=\{\,1,\ldots,n\,\}$ and $[n]'=\{\,1',\ldots,n'\,\}$, where the $\Z_2$-grading and linear ordering are induced from ${\mc N}$. 
For $a\in {\mc N}$, we write $|a|=k$ when $a$ is either $k$ or $k'$.
 
\subsection{Semistandard decomposition tableaux and Schur $P$-functions}

Let us recall the notion of semistandard decomposition tableaux \cite{GJKKK14,Se}, which is our main combinatorial object.

\begin{df}\label{SSDT}\mbox{}
\begin{itemize}
\item[(1)] A word $u=u_1 \cdots u_s$ in $\W_{\N}$ is called a {\em hook word}
if it satisfies  
$u_1 \geq u_2 \geq \cdots \geq u_k < u_{k+1} < \cdots < u_s$ for some $1 \leq k \leq s$. 
In this case, let $u\!\!\downarrow = u_1 \cdots u_k$ be the weakly decreasing subword of maximal length
and $u\!\!\uparrow = u_{k+1} \cdots u_s$ the remaining strictly increasing subword in $u$.

\item[(2)] For $\la\in \cP^+$,
let $T$ be a tableau of shifted shape $\la$ with entries in $\N$. 
Then $T$ is called a {\em semistandard decomposition tableau} of shape $\la$ if  
\begin{enumerate}
\item[(i)] $T^{(i)}$ is a hook word of length $\lambda_i$ for $1 \leq i \leq \ell(\lambda)$,  
\item[(ii)] $T^{(i)}$ is a hook subword of maximal length in $T^{(i+1)}T^{(i)}$, the concatenation of $T^{(i+1)}$ and $T^{(i)}$, for $1 \leq i < \ell(\lambda)$.
\end{enumerate}
\end{itemize}
\end{df}

For any hook word $u$, the decreasing part $u\!\!\downarrow$ is always nonempty by definition.

For $\la\in\cP^+$, let $SSDT(\la)$ be the set of semistandard decomposition tableaux of shape $\la$.
Let $x=\{x_1,x_2,\ldots\}$ be a set of formal commuting variables, and let $P_\la=P_\la(x)$ be the Schur $P$-function in $x$ corresponding to $\la\in \cP^+$ (see \cite{Mac}). 
It is shown in~\cite{Se} that $P_\la$ is given by the weight generating function of $SSDT(\la)$:
\begin{equation}\label{eq:P-function}
P_\lambda = \sum_{T\in SSDT(\la)} x^T,
\end{equation}
where $x^T = \prod_{i\geq 1} x_i^{c_i(T)}.$

\begin{rem}
Recall that the Schur $P$-function $P_\la$ can be realized as the character of tableaux $T\in SST^+_{\mc N}(\la)$ with no primed entry or entry of odd degree on the main diagonal (cf. \cite{Mac,Sa,Wo}).  
The notion of semistandard decomposition tableaux was introduced in~\cite{Se} to give a plactic monoid model for Schur $P$-functions. 
In this paper, we follow its modified version (Definition \ref{SSDT}) introduced in \cite{GJKKK14}, by which it is more easier to describe $\mf q(n)$-crystals \cite[Remark 2.6]{GJKKK14}. We also refer the reader to~\cite{CNO14} for more details on relation between the combinatorics of these two models.
\end{rem}

The following is a useful criterion for a tableau to be a semistandard decomposition one, which plays an important role in this paper.
 
\begin{prop}{\rm (\cite[Proposition 2.3]{GJKKK14})}\label{prop:SSDT-checker}
For $\la\in \cP^+$, let $T$ be a tableau of shifted shape $\la$ with entries in $\N$. 
Then $T\in SSDT(\la)$ if and only if $T^{(k)}$ is a hook word for $1\leq k\leq \ell(\la)$, 
and none of the following conditions holds for each $1\leq k < \ell(\la)$:
\begin{itemize}
\item[\rm (1)] $T_{k,1} \leq T_{k+1,i}$ for some $1 \leq i \leq \la_{k+1}$,
\item[\rm (2)] $T_{k+1,i} \geq T_{k+1,j} \geq T_{k,i+1}$ for some $1 \leq i < j \leq \la_{k+1}$,
\item[\rm (3)] $T_{k+1,j}  < T_{k,i}  < T_{k,j+1}$ for some $1 \leq i \leq j \leq \la_{k+1}$.
\end{itemize}
Equivalently, $T\in SSDT(\la)$ if and only if $T^{(k)}$ is a hook word for $1\leq k\leq \ell(\la)$, and the following conditions hold for $1\leq k < \ell(\la)$:
\begin{itemize}
\item[\rm (a)] if $T_{k,i} \leq T_{k+1,j}$ for $1 \leq i \leq j\leq  \la_{k+1}$, 
then $i\neq 1$ and $T_{k+1,i-1}\leq T_{k+1,j}$,

\item[\rm (b)] if $T_{k,i} > T_{k+1,j}$ for $1 \leq i \leq j\leq  \la_{k+1}$, 
then $T_{k,i} \geq T_{k,j+1}$.
\end{itemize}\end{prop}

For $\la\in \cP^+$, let $SSDT_n(\la)$ be the set of tableaux $T\in SSDT(\la)$ with entries in $[n]$. 
By Proposition \ref{prop:SSDT-checker}(1), we see that $SSDT_n(\la)\neq \emptyset$ if and only if $\la\in \cP^+_n$. We denote by $P_\la(x_1,\ldots,x_n)$ the Schur $P$-polynomial in $x_1,\ldots,x_n$ given by specializing $P_\la$ at $x_{n+1}=x_{n+2}=\cdots =0$. Then we have $P_\la(x_1,\ldots,x_n)=\sum_{T\in SSDT_n(\la)} x^T$.

For $\la\in\cP^+_n$,
let $H_n^\lambda$ be the element in $SSDT_n(\la)$ where the subtableau with entry  $\ell(\la)-i+1$ is a connected border strip of size $\lambda_{\ell(\lambda)-i+1}$ starting at $(i,i)\in D^+_\la$ for each $i=1,\ldots, \ell(\la)$, and 
let $L_n^\lambda$ be the one where the subtableau with entry $n-i+1$ is a connected horizontal strip of size $\lambda_{i}$ starting at $(i,i)\in D^+_\la$ for each $i=1,\ldots,\ell(\la)$.
For example, when $n=4$ and $\lambda = (4,3,1)$, we have 
\begin{displaymath}
\begin{tikzpicture}[baseline=0mm]
\def \hhh{4.4mm}
\def \vvv{5.4mm}
\draw[-] (\hhh*0,\vvv*0) rectangle (\hhh*1,\vvv*1);
\draw[-] (\hhh*1,\vvv*0) rectangle (\hhh*2,\vvv*1);
\draw[-] (\hhh*2,\vvv*0) rectangle (\hhh*3,\vvv*1);
\draw[-] (\hhh*3,\vvv*0) rectangle (\hhh*4,\vvv*1);
\draw[-] (\hhh*1,-\vvv*1) rectangle (\hhh*2,\vvv*0);
\draw[-] (\hhh*2,-\vvv*1) rectangle (\hhh*3,\vvv*0);
\draw[-] (\hhh*3,-\vvv*1) rectangle (\hhh*4,\vvv*0);
\draw[-] (\hhh*2,-\vvv*2) rectangle (\hhh*3,-\vvv*1);
\node at (\hhh*0.5,\vvv*0.5) {$3$};
\node at (\hhh*1.5,\vvv*0.5) {$2$};
\node at (\hhh*2.5,\vvv*0.5) {$2$};
\node at (\hhh*3.5,\vvv*0.5) {$1$};
\node at (\hhh*1.5,-\vvv*0.5) {$2$};
\node at (\hhh*2.5,-\vvv*0.5)  {$1$};
\node at (\hhh*3.5,-\vvv*0.5) {$1$};
\node at (\hhh*2.5,-\vvv*1.5) { $1$};
\node at (-\hhh*1.7,-\vvv*0.3) {$H_n^\lambda =$};
\end{tikzpicture}
\hskip 15mm
\begin{tikzpicture}[baseline=0mm]
\def \hhh{4.4mm}
\def \vvv{5.4mm}
\end{tikzpicture} 
\quad
\begin{tikzpicture}[baseline=0mm]
\def \hhh{4.4mm}
\def \vvv{5.4mm}
\draw[-] (\hhh*0,\vvv*0) rectangle (\hhh*1,\vvv*1);
\draw[-] (\hhh*1,\vvv*0) rectangle (\hhh*2,\vvv*1);
\draw[-] (\hhh*2,\vvv*0) rectangle (\hhh*3,\vvv*1);
\draw[-] (\hhh*3,\vvv*0) rectangle (\hhh*4,\vvv*1);
\draw[-] (\hhh*1,-\vvv*1) rectangle (\hhh*2,\vvv*0);
\draw[-] (\hhh*2,-\vvv*1) rectangle (\hhh*3,\vvv*0);
\draw[-] (\hhh*3,-\vvv*1) rectangle (\hhh*4,\vvv*0);
\draw[-] (\hhh*2,-\vvv*2) rectangle (\hhh*3,-\vvv*1);
\node at (\hhh*0.5,\vvv*0.5) {$4$};
\node at (\hhh*1.5,\vvv*0.5) {$4$};
\node at (\hhh*2.5,\vvv*0.5) {$4$};
\node at (\hhh*3.5,\vvv*0.5) {$4$};
\node at (\hhh*1.5,-\vvv*0.5) {$3$};
\node at (\hhh*2.5,-\vvv*0.5)  {$3$};
\node at (\hhh*3.5,-\vvv*0.5) {$3$};
\node at (\hhh*2.5,-\vvv*1.5) { $2$};
\node at (-\hhh*1.7,-\vvv*0.3) {$L_n^\lambda =$};
\node at (\hhh*4.5,-\vvv*2) {.};
\end{tikzpicture} 
\end{displaymath}
\noindent Indeed, $H_n^\la$ and $L_n^\la$ are the unique tableaux in $SSDT_n(\la)$ such that 
\begin{equation*}
(c_1(H^\la_n),\ldots, c_n(H^\la_n))=\la,\quad (c_1(L^\la_n),\ldots,c_n(L^\la_n))=w_0\la.
\end{equation*}
Here we assume that $\cP^+_n\subset \Z_+^n$ and the symmetric group $\mf{S}_n$ acts on $\Z_+^n$ by permutation, where $w_0$ is the longest element in $\mf S_n$.

\subsection{Crystals} 

Let us first review the crystals for the general linear Lie algebra $\gl(n)$ in \cite{Kas95, KashNaka}.

Let $P^\vee=\bigoplus_{i=1}^n\Z e_{i}$ be the dual weight lattice and $P={\rm Hom}_\Z(P^\vee,\Z)=\bigoplus_{i=1}^n\Z\epsilon_i$ the weight lattice with $\langle \epsilon_i,e_{j} \rangle =\delta_{ij}$ for $1\leq i,j\leq n$. Define a symmetric bilinear form $(\, \cdot \,|\, \cdot\,)$ on $P$ by $(\epsilon_i |\epsilon_j)=\delta_{ij}$ for $1\leq i,j\leq n$. 
Let $\{\,\alpha_i=\epsilon_i-\epsilon_{i+1}\ (i=1,\ldots,n-1)\,\}$ be the set of simple roots, and $\{\,h_i=e_{i}-e_{i+1}\ (i=1,\ldots,n-1)\,\}$ the set of simple coroots of $\gl(n)$.
Let $P^+=\{\,\la\,|\,\la\in P, \ \langle \la,h_i\rangle\geq 0 \ (i=1,\ldots,n-1) \,\}$ be the set of dominant integral weights.

A {\em $\gl(n)$-crystal} is a set
$B$ together with the maps ${\rm wt} : B \rightarrow P$,
$\varepsilon_i, \varphi_i: B \rightarrow \mathbb{Z}\cup\{-\infty\}$ and
$\te_i, \tf_i: B \rightarrow B\cup\{{\bf 0}\}$ for $i=1,\ldots,n-1$ satisfying the following conditions:  for $b\in B$ and $i=1,\ldots,n-1$,
\begin{itemize}
\item[(1)]  
$\varphi_i(b) =\langle {\rm wt}(b),h_i \rangle +
\varepsilon_i(b)$,

\item[(2)] $\varepsilon_i(\te_i b) = \varepsilon_i(b) - 1,\ \varphi_i(\te_i b) =
\varphi_i(b) + 1,\ {\rm wt}(\te_ib)={\rm wt}(b)+\alpha_i$ if $\te_i b \in B$,

\item[(3)] $\varepsilon_i(\tf_i b) = \varepsilon_i(b) + 1,\ \varphi_i(\tf_i b) =
\varphi_i(b) - 1,\ {\rm wt}({\tf_i}b)={\rm wt}(b)-\alpha_i$ if $\tf_i b \in B$,

\item[(4)] $\tf_i b = b'$ if and only if $b = \te_i b'$ for $b' \in B$,

\item[(5)] $\te_ib=\tf_ib={\bf 0}$ when $\varphi_i(b)=-\infty$.
\end{itemize}
Here ${\bf 0}$ is a formal symbol and $-\infty$ is the smallest
element in $\Z\cup\{-\infty\}$ such that $-\infty+n=-\infty$
for all $n\in\Z$. 
For $\mu\in P$, let
$B_\mu=\{\, b\in B\,  | \, {\rm wt}(b)=\mu\,\}$.
When $B_\mu$ is finite for all $\mu$, we define the character of $B$ by ${\rm ch}B=\sum_{\mu\in P}|B_\mu|e^\mu$,
where $e^\mu$ is a basis element of the group algebra $\mathbb{Q}[P]$.

Let $B_1$ and $B_2$ be $\gl(n)$-crystals.
A  tensor product $B_1\otimes B_2$
is a $\gl(n)$-crystal, which is defined to be $B_1\times B_2$  as a set with elements  denoted by
$b_1\otimes b_2$, where  

{\allowdisplaybreaks
\begin{equation*}\label{eq:tensor product of crystals}
\begin{split}
{\rm wt}(b_1\otimes b_2)&={\rm wt}(b_1)+{\rm wt}(b_2), \\
\varepsilon_i(b_1\otimes b_2)&= {\rm
max}\{\varepsilon_i(b_1),\varepsilon_i(b_2)-\langle {\rm
wt}(b_1),h_i\rangle\}, \\
\varphi_i(b_1\otimes b_2)&= {\rm max}\{\varphi_i(b_1)+\langle {\rm
wt}(b_2),h_i\rangle,\varphi_i(b_2)\},
\end{split}
\end{equation*}}
{\allowdisplaybreaks
\begin{equation}\label{eq:tensor product of crystals}
\begin{split}
{\te}_i(b_1\otimes b_2)&=
\begin{cases}
{\te}_i b_1 \otimes b_2, & \text{if $\varphi_i(b_1)\geq \varepsilon_i(b_2)$}, \\
b_1\otimes {\te}_i b_2, & \text{if
$\varphi_i(b_1)<\varepsilon_i(b_2)$},
\end{cases}\\
{\tf}_i(b_1\otimes b_2)&=
\begin{cases}
{\tf}_i b_1 \otimes b_2, & \text{if  $\varphi_i(b_1)>\varepsilon_i(b_2)$}, \\
b_1\otimes {\tf}_i b_2, & \text{if $\varphi_i(b_1)\leq
\varepsilon_i(b_2)$},
\end{cases}
\end{split}
\end{equation}}
\noindent for $i=1,\ldots,n-1$. Here we assume that ${\bf 0}\otimes
b_2=b_1\otimes {\bf 0}={\bf 0}$.

For $\la\in \cP_n$, let $B_n(\la)$ be the crystal associated to an irreducible $\gl(n)$-module with highest weight $\la$, where we regard $\la$ as $\sum_{i=1}^n\la_i\epsilon_i\in P^+$.
We may regard $[n]$ as $B_n(\epsilon_1)$, where ${\rm wt}(k)=\epsilon_k$ for $k\in[n]$, and hence $\W_{[n]}$ as a $\gl(n)$-crystal
where we identify $w=w_1\ldots w_r$ with $w_1\otimes \cdots \otimes w_r\in B_n(\epsilon_1)^{\otimes r}$. 
The crystal structure on $\W_{[n]}$ is easily described by so-called the {signature rule} (cf. \cite[Section 2.1]{KashNaka}).
For $\la\in \cP_n$, the set $SST_{[n]}(\la)$ becomes a $\gl(n)$-crystal under the identification of $T$ with $w(T)\in \W_{[n]}$, and it is isomorphic to $B_n(\la)$ \cite{KashNaka}. In general, one can define a $\gl(n)$-crystal structure on $SST_{[n]}(\la/\mu)$ for a skew diagram $\la/\mu$. By abuse of notation, we set $B_n(\la/\mu):=SST_{[n]}(\la/\mu)$.
\vskip 2mm

Next, let us review the notion of crystals associated to polynomial representations of the queer Lie superalgebra ${\mf q}(n)$ developed in \cite{GJKKK14,GJKKK15}.  
\begin{df}\label{def:q(n)-crystals}
A {\em $\mathfrak{q}(n)$-crystal} is a set $B$ together with the maps ${\rm wt} : B \rightarrow P$,
$\varepsilon_i, \varphi_i: B \rightarrow \mathbb{Z}\cup\{-\infty\}$ and
$\te_i, \tf_i: B \rightarrow B\cup\{{\bf 0}\}$ for $i\in I:=\{\,1,\ldots,n-1, \ov{1}\,\}$ satisfying the following conditions: 
\begin{itemize}
\item[(1)] $B$ is a $\gl(n)$-crystal with respect to ${\rm wt}$, $\varepsilon_i$, $\varphi_i$, $\te_i, \tf_i$ for $i=1,\ldots, n-1$,

\item[(2)] $\wt(b) \in \bigoplus_{i\in [n]}\Z_+\epsilon_i$ for $b\in B$, 

\item[(3)] $\wt(\widetilde{e}_{\overline{1}}b) = \wt(b) + \alpha_1$, $\wt(\widetilde{f}_{\overline{1}}b)=\wt(b)-\alpha_1$ for $b\in B$,

\item[(4)] $\widetilde{f}_{\overline{1}}b = b'$ if and only if 
$b = \widetilde{e}_{\overline{1}}b'$ for all $b, b' \in B$,

\item[(5)] for $3 \leq i \leq n-1$, we have
\begin{itemize}
\item[(i)] the operators $\widetilde{e}_{\overline{1}}$ and $\widetilde{f}_{\overline{1}}$ commute with $\widetilde{e}_{i}$, $\widetilde{f}_{i}$,

\item[(ii)] if $\widetilde{e}_{\overline{1}}b \in B$, then $\varepsilon_i(\widetilde{e}_{\overline{1}}b)=\varepsilon_i(b)$ and 
$\varphi_i(\widetilde{e}_{\overline{1}}b)=\varphi_i(b)$.
\end{itemize}
\end{itemize}
\end{df}

Let $\B_n$ be a $\q(n)$-crystal which is the $\gl(n)$-crystal $B_n(\epsilon_1)$ together with $\tf_{\ov 1} \ \boxed{1}=\boxed{2}$ (in dashed arrow):
\begin{center}
\begin{tikzpicture}[baseline=0mm]
\def \hhh{4.4mm}
\def \hhhh{18mm}
\def \cdotssize{5mm}
\def \vvv{5mm}
\draw[-] (\hhh*0,\vvv*0) rectangle (\hhh*1,\vvv*1);
\draw[-] (\hhh*0+\hhhh*1,\vvv*0) rectangle (\hhh*1+\hhhh*1,\vvv*1);
\draw[-] (\hhh*0+\hhhh*2,\vvv*0) rectangle (\hhh*1+\hhhh*2,\vvv*1);
\draw[-] (\cdotssize+\hhh*0+\hhhh*4,\vvv*0) rectangle (\cdotssize+\hhh*1+\hhhh*4,\vvv*1);
\node at (\hhh*0.5,\vvv*0.5) {$1$};
\node at (\hhh*0.5+\hhhh,\vvv*0.5) {$2$};
\node at (\hhh*0.5+\hhhh*2,\vvv*0.5) {$3$};
\node at (\cdotssize*0.5+\hhh*0.5+\hhhh*3,\vvv*0.5) {$\cdots$};
\node at (\cdotssize+\hhh*0.5+\hhhh*4,\vvv*0.5) {$n$};
\draw[->] (\hhh*1.3,\vvv*0.7) -- (-\hhh*0.3+\hhhh*1,\vvv*0.7) node[midway,above] () {\tiny $1$};
\draw[->,densely dashed] (\hhh*1.3,\vvv*0.3) -- (-\hhh*0.3+\hhhh*1,\vvv*0.3) node[midway,below] () {\tiny $\overline{1}$};
\draw[->] (\hhh*1.3+\hhhh*1,\vvv*0.5) -- (-\hhh*0.3+\hhhh*2,\vvv*0.5) node[midway,above] () {\tiny $2$};
\draw[->] (\hhh*1.3+\hhhh*2,\vvv*0.5) -- (-\hhh*0.3+\hhhh*3,\vvv*0.5) node[midway,above] () {\tiny $3$};
\draw[->] (\cdotssize+\hhh*1.1+\hhhh*3,\vvv*0.5) -- (\cdotssize-\hhh*0.3+\hhhh*4,\vvv*0.5) node[midway,above] () {\tiny $n-1$};
\node at (\cdotssize+\hhh*1.3+\hhhh*4,\vvv*0.1) {.};
\end{tikzpicture}
\end{center}
Here we write $b\,\stackrel{i}{\longrightarrow}\,b'$ if $\tf_i b= b'$ for $b,b'\in B$ and $i\in I\setminus\{\ov{1}\}$ as usual, and $b\ \stackrel{\ov 1}{\dashrightarrow}\ b'$ if $\tf_{\ov{1}} b= b'$.

For $\mathfrak{q}(n)$-crystals $B_1$ and $B_2$,
the tensor product $B_1 \otimes B_2$ is the $\gl(n)$-crystal $B_1 \otimes B_2$
where the actions of
$\widetilde{e}_{\overline{1}}$ and $\widetilde{f}_{\overline{1}}$
are given by 
\begin{equation}\label{eq:tensor product of qn-crystals}
\begin{split}
{\te}_{\ov1}(b_1\otimes b_2)&=
\begin{cases}
{\te}_{\ov1}b_1 \otimes b_2, & \text{if } \langle \epsilon_1,\wt(b_2)\rangle = \langle \epsilon_2,\wt(b_2)\rangle=0,\\
b_1\otimes {\te}_{\ov1}b_2, & \text{otherwise}, 
\end{cases}\\
{\tf}_{\ov1}(b_1\otimes b_2)&=
\begin{cases}
{\tf}_{\ov1} b_1 \otimes b_2, & \text{if }  \langle \epsilon_1,\wt(b_2)\rangle = \langle \epsilon_2,\wt(b_2)\rangle=0, \\
b_1\otimes {\tf}_{\ov1}b_2, & \text{otherwise.}
\end{cases}
\end{split}
\end{equation}
\noindent Then it is easy to see that $B_1\otimes B_2$ is a $\mathfrak{q}(n)$-crystal. In particular, $\W_{[n]}$ is also a $\q(n)$-crystal.

Let $B$ be a $\q(n)$-crystal. Suppose that $B$ is a regular $\gl(n)$-crystal, that is, each connected component in $B$ is isomorphic to $B_n(\la)$ for some $\la\in\cP_n$. Let $W=\mf S_n$ be the Weyl group of $\gl(n)$ which is generated by the simple reflection $r_i$ corresponding to $\alpha_i$ for $i=1,\ldots,n-1$. We have a group action of $W$ on $B$ 
denoted by $S$ such that 
\begin{equation*}
S_{r_i}(b)=
\begin{cases}
\tf_i^{\langle {\rm wt}(b),h_i\rangle}b,& \text{if $\langle {\rm wt}(b),h_i\rangle\geq 0$},\\
\te_i^{\,-\langle {\rm wt}(b),h_i\rangle}b,& \text{if $\langle {\rm wt}(b),h_i\rangle\leq 0$},
\end{cases}
\end{equation*}
for $b\in B$ and $i=1,\ldots,n-1$. 
For $2\leq i\leq n-1$, let $w_i\in W$ be such that $w_i(\alpha_{i})=\alpha_1$, and 
let  
\begin{equation}
\te_{\ov{i}}=S_{w_i^{-1}}\te_{\ov{1}}S_{w_i},\quad
\tf_{\ov{i}}=S_{w_i^{-1}}\tf_{\ov{1}}S_{w_i}.
\end{equation}
For $b\in B$, we say that $b$ is a {\em $\q(n)$-highest weight vector} if 
$\te_i b=\te_{\ov i}b={\bf 0}$ for $1\leq i \leq n-1$, and 
$b$ is a {\em $\q(n)$-lowest weight vector} if $S_{w_0}b$ is a $\q(n)$-highest weight vector.

For $\la\in \cP^+$, let $\B_n(\la)=SSDT_n(\la)$, and consider an injective map
\begin{equation}\label{eq:embedding of SSDT}
\xymatrixcolsep{2.5pc}\xymatrixrowsep{0pc}
\xymatrix{
\B_n(\la)\ \  \ar@{^{(}->}[r]  & \ \W_{[n]} \\
\ \ \ T \ \  \ \   \ar@{|->}[r]   &  \  \makebox[2em][l]{$w_\rev(T)$.}}
\end{equation}
Then we have the following.

\begin{thm}{\rm (\cite[Theorem 2.5]{GJKKK14})}
Let $\lambda \in \cP^+_n$ be given.
\begin{itemize}
\item[\rm (a)] The image of $\B_n(\lambda)$ in \eqref{eq:embedding of SSDT} together with $\{{\bf 0}\}$ is invariant under the action of $\te_i$ and $\tf_i$ for $i\in I$, and hence $\B_n(\la)$ is a $\q(n)$-crystal.

\item[\rm (b)] The $\q(n)$-crystal $\B_n(\la)$ is connected where $H^\lambda_n$ is a unique $\q(n)$-highest weight vector and $L^\lambda_n$ is a unique $\q(n)$-lowest weight vector.
\end{itemize}
\end{thm}

\begin{rem}
In \cite{GJKKK15}, a semisimple tensor category over the quantum superalgebra $U_q(\q(n))$ is introduced, and it is shown that each irreducible highest weight module $V_n(\la)$ in this category, parametrized by $\la\in\cP^+_n$, has a crystal base. Furthermore, it is shown in \cite[Theorem 2.5(c)]{GJKKK14} that the crystal of $V_n(\la)$ is isomorphic to $\B_n(\la)$.
\end{rem}

\begin{figure}[b]
\begin{tikzpicture}[scale=1,xscale=1,yscale=0.7]
\def \hhh{70bp}
\def \vvv{75bp}
\def \widthh{.6ex}
\def \highhh{-.3ex}
\node (node_13) at (0bp,0bp) [draw,draw=none] {${\def\lr#1{\multicolumn{1}{|@{\hspace{\widthh}}c@{\hspace{\widthh}}|}{\raisebox{\highhh}{$#1$}}}\raisebox{\highhh}
{$\begin{array}[b]{*{3}c}\cline{1-3}\lr{2}&\lr{1}&\lr{1}\\\cline{1-3}&\lr{1}\\\cline{2-2}\end{array}$}}$};
\node (node_22) at (-\hhh,-\vvv) [draw,draw=none] {${\def\lr#1{\multicolumn{1}{|@{\hspace{\widthh}}c@{\hspace{\widthh}}|}{\raisebox{\highhh}{$#1$}}}\raisebox{\highhh}
{$\begin{array}[b]{*{3}c}\cline{1-3}\lr{2}&\lr{2}&\lr{1}\\\cline{1-3}&\lr{1}\\\cline{2-2}\end{array}$}}$};
\node (node_23) at (0,-\vvv) [draw,draw=none] {${\def\lr#1{\multicolumn{1}{|@{\hspace{\widthh}}c@{\hspace{\widthh}}|}{\raisebox{\highhh}{$#1$}}}\raisebox{\highhh}
{$\begin{array}[b]{*{3}c}\cline{1-3}\lr{3}&\lr{1}&\lr{1}\\\cline{1-3}&\lr{1}\\\cline{2-2}\end{array}$}}$};
\node (node_24) at (0+\hhh,-\vvv) [draw,draw=none] {${\def\lr#1{\multicolumn{1}{|@{\hspace{\widthh}}c@{\hspace{\widthh}}|}{\raisebox{\highhh}{$#1$}}}\raisebox{\highhh}
{$\begin{array}[b]{*{3}c}\cline{1-3}\lr{2}&\lr{1}&\lr{2}\\\cline{1-3}&\lr{1}\\\cline{2-2}\end{array}$}}$};
\node (node_31) at (0-\hhh*2,-\vvv*2) [draw,draw=none] {${\def\lr#1{\multicolumn{1}{|@{\hspace{\widthh}}c@{\hspace{\widthh}}|}{\raisebox{\highhh}{$#1$}}}\raisebox{\highhh}
{$\begin{array}[b]{*{3}c}\cline{1-3}\lr{2}&\lr{2}&\lr{2}\\\cline{1-3}&\lr{1}\\\cline{2-2}\end{array}$}}$};
\node (node_32) at (0-\hhh,-\vvv*2) [draw,draw=none] {${\def\lr#1{\multicolumn{1}{|@{\hspace{\widthh}}c@{\hspace{\widthh}}|}{\raisebox{\highhh}{$#1$}}}\raisebox{\highhh}{$\begin{array}[b]{*{3}c}\cline{1-3}\lr{3}&\lr{2}&\lr{1}\\\cline{1-3}&\lr{1}\\\cline{2-2}\end{array}$}}$};
\node (node_33) at (0,-\vvv*2) [draw,draw=none] {${\def\lr#1{\multicolumn{1}{|@{\hspace{\widthh}}c@{\hspace{\widthh}}|}{\raisebox{\highhh}{$#1$}}}\raisebox{\highhh}
{$\begin{array}[b]{*{3}c}\cline{1-3}\lr{3}&\lr{1}&\lr{1}\\\cline{1-3}&\lr{2}\\\cline{2-2}\end{array}$}}$};
\node (node_34) at (\hhh,-\vvv*2) [draw,draw=none] {${\def\lr#1{\multicolumn{1}{|@{\hspace{\widthh}}c@{\hspace{\widthh}}|}{\raisebox{\highhh}{$#1$}}}\raisebox{\highhh}
{$\begin{array}[b]{*{3}c}\cline{1-3}\lr{3}&\lr{1}&\lr{2}\\\cline{1-3}&\lr{1}\\\cline{2-2}\end{array}$}}$};
\node (node_36) at (\hhh*3,-\vvv*2) [draw,draw=none] {${\def\lr#1{\multicolumn{1}{|@{\hspace{\widthh}}c@{\hspace{\widthh}}|}{\raisebox{\highhh}{$#1$}}}\raisebox{\highhh}
{$\begin{array}[b]{*{3}c}\cline{1-3}\lr{2}&\lr{1}&\lr{3}\\\cline{1-3}&\lr{1}\\\cline{2-2}\end{array}$}}$};
\node (node_41) at (-\hhh*2,-\vvv*3) [draw,draw=none] {${\def\lr#1{\multicolumn{1}{|@{\hspace{\widthh}}c@{\hspace{\widthh}}|}{\raisebox{\highhh}{$#1$}}}\raisebox{\highhh}
{$\begin{array}[b]{*{3}c}\cline{1-3}\lr{3}&\lr{2}&\lr{2}\\\cline{1-3}&\lr{1}\\\cline{2-2}\end{array}$}}$};
\node (node_42) at (-\hhh,-\vvv*3) [draw,draw=none] {${\def\lr#1{\multicolumn{1}{|@{\hspace{\widthh}}c@{\hspace{\widthh}}|}{\raisebox{\highhh}{$#1$}}}\raisebox{\highhh}
{$\begin{array}[b]{*{3}c}\cline{1-3}\lr{3}&\lr{3}&\lr{1}\\\cline{1-3}&\lr{1}\\\cline{2-2}\end{array}$}}$};
\node (node_43) at (0,-\vvv*3) [draw,draw=none] {${\def\lr#1{\multicolumn{1}{|@{\hspace{\widthh}}c@{\hspace{\widthh}}|}{\raisebox{\highhh}{$#1$}}}\raisebox{\highhh}
{$\begin{array}[b]{*{3}c}\cline{1-3}\lr{3}&\lr{2}&\lr{1}\\\cline{1-3}&\lr{2}\\\cline{2-2}\end{array}$}}$};
\node (node_44) at (\hhh,-\vvv*3) [draw,draw=none] {${\def\lr#1{\multicolumn{1}{|@{\hspace{\widthh}}c@{\hspace{\widthh}}|}{\raisebox{\highhh}{$#1$}}}\raisebox{\highhh}
{$\begin{array}[b]{*{3}c}\cline{1-3}\lr{3}&\lr{1}&\lr{2}\\\cline{1-3}&\lr{2}\\\cline{2-2}\end{array}$}}$};
\node (node_45) at (\hhh*2,-\vvv*3) [draw,draw=none] {${\def\lr#1{\multicolumn{1}{|@{\hspace{\widthh}}c@{\hspace{\widthh}}|}{\raisebox{\highhh}{$#1$}}}\raisebox{\highhh}
{$\begin{array}[b]{*{3}c}\cline{1-3}\lr{3}&\lr{1}&\lr{3}\\\cline{1-3}&\lr{1}\\\cline{2-2}\end{array}$}}$};
\node (node_46) at (\hhh*3,-\vvv*3) [draw,draw=none] {${\def\lr#1{\multicolumn{1}{|@{\hspace{\widthh}}c@{\hspace{\widthh}}|}{\raisebox{\highhh}{$#1$}}}\raisebox{\highhh}
{$\begin{array}[b]{*{3}c}\cline{1-3}\lr{2}&\lr{2}&\lr{3}\\\cline{1-3}&\lr{1}\\\cline{2-2}\end{array}$}}$};
\node (node_51) at (-\hhh*2,-\vvv*4) [draw,draw=none] {${\def\lr#1{\multicolumn{1}{|@{\hspace{\widthh}}c@{\hspace{\widthh}}|}{\raisebox{\highhh}{$#1$}}}\raisebox{\highhh}
{$\begin{array}[b]{*{3}c}\cline{1-3}\lr{3}&\lr{3}&\lr{2}\\\cline{1-3}&\lr{1}\\\cline{2-2}\end{array}$}}$};
\node (node_52) at (-\hhh,-\vvv*4) [draw,draw=none] {${\def\lr#1{\multicolumn{1}{|@{\hspace{\widthh}}c@{\hspace{\widthh}}|}{\raisebox{\highhh}{$#1$}}}\raisebox{\highhh}
{$\begin{array}[b]{*{3}c}\cline{1-3}\lr{3}&\lr{3}&\lr{1}\\\cline{1-3}&\lr{2}\\\cline{2-2}\end{array}$}}$};
\node (node_53) at (0,-\vvv*4) [draw,draw=none] {${\def\lr#1{\multicolumn{1}{|@{\hspace{\widthh}}c@{\hspace{\widthh}}|}{\raisebox{\highhh}{$#1$}}}\raisebox{\highhh}
{$\begin{array}[b]{*{3}c}\cline{1-3}\lr{3}&\lr{2}&\lr{2}\\\cline{1-3}&\lr{2}\\\cline{2-2}\end{array}$}}$};
\node (node_54) at (\hhh*1,-\vvv*4) [draw,draw=none] {${\def\lr#1{\multicolumn{1}{|@{\hspace{\widthh}}c@{\hspace{\widthh}}|}{\raisebox{\highhh}{$#1$}}}\raisebox{\highhh}
{$\begin{array}[b]{*{3}c}\cline{1-3}\lr{3}&\lr{1}&\lr{3}\\\cline{1-3}&\lr{2}\\\cline{2-2}\end{array}$}}$};
\node (node_56) at (\hhh*3,-\vvv*4) [draw,draw=none] {${\def\lr#1{\multicolumn{1}{|@{\hspace{\widthh}}c@{\hspace{\widthh}}|}{\raisebox{\highhh}{$#1$}}}\raisebox{\highhh}
{$\begin{array}[b]{*{3}c}\cline{1-3}\lr{3}&\lr{2}&\lr{3}\\\cline{1-3}&\lr{1}\\\cline{2-2}\end{array}$}}$};
\node (node_62) at (-\hhh,-\vvv*5) [draw,draw=none] {${\def\lr#1{\multicolumn{1}{|@{\hspace{\widthh}}c@{\hspace{\widthh}}|}{\raisebox{\highhh}{$#1$}}}\raisebox{\highhh}{$\begin{array}[b]{*{3}c}\cline{1-3}\lr{3}&\lr{3}&\lr{3}\\\cline{1-3}&\lr{1}\\\cline{2-2}\end{array}$}}$};
\node (node_63) at (0,-\vvv*5) [draw,draw=none] {${\def\lr#1{\multicolumn{1}{|@{\hspace{\widthh}}c@{\hspace{\widthh}}|}{\raisebox{\highhh}{$#1$}}}\raisebox{\highhh}
{$\begin{array}[b]{*{3}c}\cline{1-3}\lr{3}&\lr{3}&\lr{2}\\\cline{1-3}&\lr{2}\\\cline{2-2}\end{array}$}}$};
\node (node_64) at (\hhh,-\vvv*5) [draw,draw=none] {${\def\lr#1{\multicolumn{1}{|@{\hspace{\widthh}}c@{\hspace{\widthh}}|}{\raisebox{\highhh}{$#1$}}}\raisebox{\highhh}
{$\begin{array}[b]{*{3}c}\cline{1-3}\lr{3}&\lr{2}&\lr{3}\\\cline{1-3}&\lr{2}\\\cline{2-2}\end{array}$}}$};
\node (node_73) at (0,-\vvv*6) [draw,draw=none] {${\def\lr#1{\multicolumn{1}{|@{\hspace{\widthh}}c@{\hspace{\widthh}}|}{\raisebox{\highhh}{$#1$}}}\raisebox{\highhh}
{$\begin{array}[b]{*{3}c}\cline{1-3}\lr{3}&\lr{3}&\lr{3}\\\cline{1-3}&\lr{2}\\\cline{2-2}\end{array}$}}$};
\definecolor{strokecol}{rgb}{0.0,0.0,0.0};
\pgfsetstrokecolor{strokecol}
\draw [red,->] (node_13) -- (node_22)  node[midway,above,sloped] {\tiny $1$};
\draw [blue,->] (node_13) -- node[right] {\tiny $2$} ++ (node_23);
\draw [purple,->,densely dashed] (node_13) -- node[above,midway,sloped] {\tiny $\overline{1}$} ++ (node_24);
\draw [red,->,transform canvas={xshift=-3pt,yshift=4pt}] (node_22) -- node[midway,above,sloped] {\tiny $1$} ++ (node_31);
\draw [red,->,densely dashed] (node_22) -- (node_31) node[midway,below,sloped] {\tiny $\overline{1}$};
\draw [blue,->] (node_22) -- node[right] {\tiny $2$} ++ (node_32);
\draw [red,->] (node_23) -- node[right] {\tiny $1$} ++ (node_33);
\draw [purple,densely dashed,->] (node_23) -- node[midway,above,sloped] {\tiny $\overline{1}$} ++ (node_34);
\draw [blue,->] (node_24) -- node[right] {\tiny $2$} ++ (node_34);
\draw [blue,->] (node_31) -- node[left] {\tiny $2$} ++ (node_41);
\draw [red,->,transform canvas={xshift=-3pt,yshift=4pt}] (node_32) -- node[above,sloped] {\tiny $1$} ++ (node_41);
\draw [purple,->,densely dashed] (node_32) -- node[below,sloped] {\tiny $\ov 1$} ++ (node_41);
\draw [blue,->] (node_32) -- node[right] {\tiny $2$} ++ (node_42);
\draw [red,->] (node_33) -- node[right] {\tiny $1$} ++ (node_43);
\draw [purple,densely dashed,->] (node_33) -- node[midway,above,sloped] {\tiny $\overline{1}$} ++ (node_44);
\draw [red,->] (node_34) -- node[right] {\tiny $1$} ++ (node_44);
\draw [blue,->] (node_34) -- node[midway,above,sloped] {\tiny $2$} ++ (node_45);
\draw [red,->,transform canvas={xshift=-3pt}] (node_36) -- node[left] {\tiny $1$} ++ (node_46);
\draw [purple,->,densely dashed,transform canvas={xshift=3pt}] (node_36) -- node[right] {\tiny $\overline{1}$} ++ (node_46);
\draw [blue,->] (node_41) -- node[left] {\tiny $2$} ++  (node_51);
\draw [red,->] (node_42) -- node[left] {\tiny $1$} ++ (node_52);
\draw [blue,->] (node_43) -- node[midway,above,sloped] {\tiny $2$} ++ (node_52);
\draw [red,->,transform canvas={xshift=-3pt}] (node_43) -- node[left] {\tiny $1$} ++ (node_53);
\draw [purple,->,densely dashed,transform canvas={xshift=3pt}] (node_43) -- node[right] {\tiny $\ov 1$} ++ (node_53);
\draw [purple,densely dashed,->] (node_42) -- node[above,sloped] {\tiny $\overline{1}$} ++ (node_51);
\draw [blue,->] (node_44) -- node[right] {\tiny $2$} ++  (node_54);
\draw [red,->] (node_45) -- node[above,sloped] {\tiny $1$} ++ (node_54);
\draw [purple,densely dashed,->] (node_45) -- node[above,sloped] {\tiny $\overline{1}$} ++ (node_56);
\draw [blue,->] (node_46) -- node[right] {\tiny $2$} ++  (node_56);
\draw [blue,->] (node_51) -- node[midway,above,sloped] {\tiny $2$} ++ (node_62);
\draw [red,->,transform canvas={xshift=3pt,yshift=4pt}] (node_52) -- node[midway,above,sloped] {\tiny $1$} ++ (node_63);
\draw [purple,->,densely dashed] (node_52) -- node[midway,below,sloped] {\tiny $1$} ++ (node_63);
\draw [blue,->] (node_53) -- node[right] {\tiny $2$} ++ (node_63);
\draw [red,->,transform canvas={xshift=-3pt}] (node_54) -- node[left] {\tiny $1$} ++ (node_64);
\draw [purple,->,densely dashed,transform canvas={xshift=3pt}] (node_54) -- node[right] {\tiny $\overline{1}$} ++ (node_64);
\draw [red,->,transform canvas={xshift=3pt,yshift=4pt}] (node_62) -- node[above,midway,sloped] {\tiny $1$} ++ (node_73);
\draw [purple,->,densely dashed] (node_62) -- node[below,midway,sloped] {\tiny $\overline{1}$} ++ (node_73);
\draw [blue,->] (node_63) -- node[right] {\tiny $2$} ++ (node_73);
\end{tikzpicture}
\caption{The $\mathfrak{q}(3)$-crystal $\B_3(3,1)$}
\label{Fig:q(3)-crystal}
\end{figure}
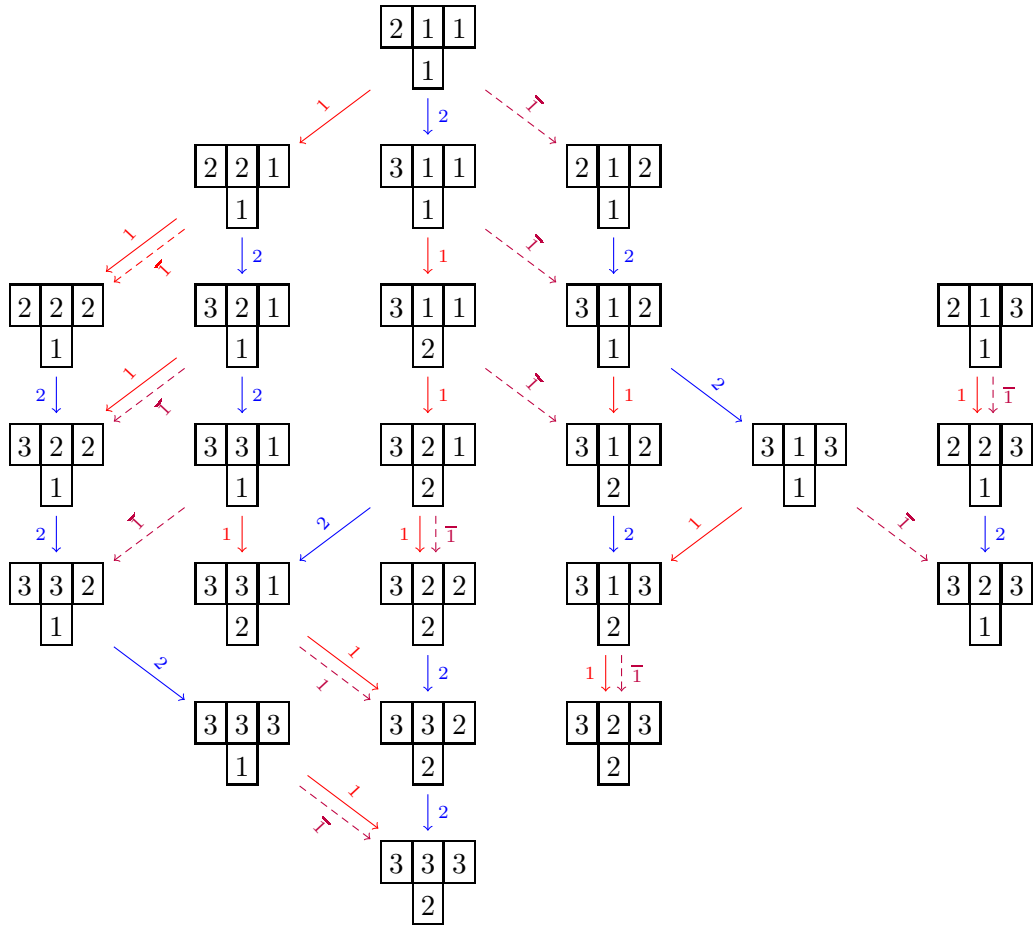

Let $B_1$ and $B_2$ be $\q(n)$-crystals. For $b_1\in B_1$ and $b_2\in B_2$, let us say that $b_1$ and $b_2$ are equivalent and write $b_1\equiv b_2$ if there exists an isomorphism of $\q(n)$-crystals $\psi : C(b_1) \longrightarrow C(b_2)$ such that $\psi(b_1)=b_2$ where $C(b_i)$ denotes the connected component of $b_i\in B_i$ $(i=1,2)$ as a $\q(n)$-crystal. 

By \cite[Theorem 4.6]{GJKKK15}, each connected component in $\B_n^{\otimes N}$ $(N\geq 1)$ is isomorphic to $\B_n(\la)$ for some $\la\in \cP^+_n$ with $|\la|=N$.
Indeed, for $b=b_1 \otimes   \cdots \otimes b_N\in \B_n^{\otimes N}$, there exists a unique $\la\in\cP_n^+$ and $T\in \B_n(\la)$ such that $b\equiv T$. In particular, $b$ is a $\q(n)$-lowest (resp. $\q(n)$-highest) weight vector if and only if $b\equiv L^\la_n$ (resp. $H^\la_n$).

The following lemma plays a crucial role in characterization of $\q(n)$-lowest weight vectors in $\B_n^{\otimes N}$ and hence describing the decompositions of $\B_n^{\otimes N}$ and $\B_n(\mu)\otimes\B_n(\nu)$ ($\mu,\nu\in\cP_n^+$) into connected components in \cite{GJKKK14}.

\begin{lem}{\rm (\cite[Lemma 1.15, Corollary 1.16]{GJKKK14})}
\label{lem:q(n)-lattice property}
For $b=b_1 \otimes \cdots \otimes b_N\in \B_n^{\otimes N}$, the following are equivalent:
\begin{itemize}
\item[\rm (1)] $b$ is a $\q(n)$-lowest weight vector,

\item[\rm (2)] $b'=b_2 \otimes \cdots \otimes b_N$ is a $\q(n)$-lowest weight vector and $\epsilon_{b_1}+{\rm wt}(b')\in w_0\cP_n^+$,

\item[\rm (3)] $\wt(b_M \otimes \cdots \otimes b_N) \in w_0 \cP^+_n$ for all $1\leq M\leq N$.
\end{itemize}
\end{lem} 

Hence, we have the following immediately by Lemma \ref{lem:q(n)-lattice property}.

\begin{cor}\label{cor:lattice property in general}
For $\la^{(1)},\ldots,\la^{(s)}\in\cP_n^+$ and $T_1 \otimes \cdots \otimes T_s\in\B_n(\la^{(1)})\otimes \cdots \otimes \B_n(\la^{(s)})$, the following are equivalent: 
\begin{itemize}
\item[\rm (1)] $T_1 \otimes \cdots \otimes T_s$ is a $\q(n)$-lowest weight vector, 

\item[\rm (2)] $T_r\otimes \cdots\otimes T_s\in \B_n(\la^{(s)})\otimes \cdots \otimes \B_n(\la^{(r)})$ is a $\q(n)$-lowest weight vector for all $1\leq r\leq s$.
\end{itemize}
\end{cor}

Note that we do not have an analogue of Lemma \ref{lem:q(n)-lattice property} for $\q(n)$-highest weight vectors.

\begin{rem}\label{rem:stability}
Let $m\geq n$ be a positive integer, and put $t=m-n$. For $N\geq 1$, let $\psi_t : \B_n^{\otimes N} \longrightarrow \B_m^{\otimes N}$ be the map given by 
$\psi_t(u_1\otimes \cdots \otimes u_N)=(u_1+t)\otimes \cdots \otimes (u_N+t)$. 
Then  for $\la\in\cP_n^+$ and $u\in \B_n^{\otimes N}$ we have $u\equiv L_n^\la$ if and only if $\psi_t(u)\equiv L_m^\la$. This implies that the multiplicity of $\B_n(\la)$ in $\B_n^{\otimes N}$ is equal to that of $\B_m(\la)$ in $\B_m^{\otimes N}$ for $\la\in \cP_n^+$.
\end{rem}

\section{Littlewood-Richardson rule for Schur $P$-functions}

For $\la, \mu, \nu \in \cP^+$, the {\em shifted Littlewood-Richardson coefficients} $f^{\la}_{ \mu \nu}$ are the coefficients given by
\begin{equation}\label{eq:LRS for P-functions}
P_{\mu}  P_{\nu} 
=\underset{\la}{\sum}f^{\la}_{\mu \nu}P_{\la}.
\end{equation}
In this section we give a new combinatorial description of $f^\la_{\mu\nu}$ using the theory of $\q(n)$-crystals. We also show that our description of $f^\la_{\mu\nu}$ is equivalent to the Stembridge's description \cite{Ste}.

\subsection{Shifted Littlewood-Richardson rule}

\begin{df}\label{Def:w and w^*}
{\rm
Let $w=w_1\cdots w_N$ be a word in $\W_{\mc N}$. Let $m_k=c_k(w)+c_{k'}(w)$ for $k\geq 1$. We define $w^*=w^*_1\cdots w^*_N$ to be the word obtained from $w$ after applying the following steps for each $k\geq 1$:

\begin{itemize}

\item[(1)] Consider the letters $w_i$'s with $|w_i|=k$. Label them with $1,2,\ldots, m_k$ (as subscripts), first enumerating the $w_p$'s  with $w_p=k$ from left to right, and then $w_q$'s with $w_q=k'$ from right to left.

\item[(2)] After the step (1), remove all \ $'$\ \ in each labeled letter $k'_j$, that is, replace any $k'_j$ with $k_j$ for $c_{k}(w)< j\leq m_k$.

\end{itemize}
}
\end{df}
\begin{exm}{\rm
\begin{equation*}
w=11'11'1  \quad\longrightarrow \quad 1_11'_51_21'_41_3 \quad\longrightarrow \quad  w^*=1_11_51_21_41_3
\end{equation*}
\begin{equation*}
w=21'12'2'121 \quad\longrightarrow \quad 
2_{\red 1}1'_{\blue 4}1_{\blue 1}2'_{\red 4}2'_{\red 3}1_{\blue 2}2_{\red 2}1_{\blue 3}
\quad\longrightarrow \quad w^*=2_{\red 1}1_{\blue 4}1_{\blue 1}2_{\red 4}2_{\red 3}1_{\blue 2}2_{\red 2}1_{\blue 3}
\end{equation*}
}
\end{exm}

\begin{df}
\label{Def:new lattice property}
{\rm
Let $w=w_1\cdots w_N\in \W_{\mc N}$ be given.
We say that $w$ satisfies the {\em ``lattice property"} if 
the word $w^*=w^*_1\cdots w^*_N$ associated to $w$ given in Definition \ref{Def:w and w^*} satisfies the following for $k\geq 1$:
\begin{itemize}
\item[(L1)] if $w^*_i=k_1$, then no $k+1_j$ for $j\geq 1$ occurs in $w^*_1\cdots w^*_{i-1}$,

\item[(L2)] if $(w^*_s,w^*_t)=(k+1_i,k_{i+1})$ for some $s<t$ and $i \geq 1$, 
then no $k+1_j$ for $i<j$ occurs in $w^*_s\cdots w^*_t$, 

\item[(L3)] if $(w^*_s,w^*_t)=(k_{j+1},k+1_{j})$ for some $s<t$ and $j\geq 1$, 
then no $k_i$ for $i\leq j$ occurs in $w^*_s\cdots w^*_t$.
\end{itemize}}
\end{df}

\begin{df}\label{Def:new LR tab}
{\rm
For $\la, \mu, \nu \in \cP^+$,
let $\tt{F}^\la_{\mu\nu}$ be the set of tableaux $Q$ such that  
\begin{itemize}
\item[(1)] $Q\in SST^+_{\mc N}(\la/\mu)$ with $c_{k}(Q)+c_{k'}(Q)=\nu_k$ for $k\geq 1$,

\item[(2)] for $k\geq 1$, if $x$ is the rightmost letter in $w(Q)$ with $|x|=k$, then $x=k$,

\item[(3)] $w(Q)$ satisfies the ``lattice property" in Definition \ref{Def:new lattice property}. 

\end{itemize}}
\end{df}
Then we have the following characterization of $f^\la_{\mu\nu}$.

\begin{thm}\label{Thm:our-characterization}
For $\la, \mu, \nu \in \cP^+$, we have 
$$
f^\la_{\mu\nu} = \left|\tt{F}^\la_{\mu\nu}\right|,
$$
that is, the shifted LR coefficient $f^\la_{\mu \nu}$ is equal to the number of tableaux in $\tt{F}^\la_{\mu\nu}$.
\end{thm}
\pf Choose $n$ such that $\la, \mu, \nu\in \cP_n^+$. Put
\begin{equation}\label{eq:lowest weight vectors for LRS}
{\tt L}^\la_{\mu\nu}
=\{\,T\,|\, T \in \B_n(\nu),\ T \otimes L^\mu_n \equiv L^\la_n\,\}.
\end{equation}
By Corollary \ref{cor:lattice property in general},
we have
\begin{equation}\label{eq:tensor decomp of mu nu}
\B_n(\nu)\otimes \B_n(\mu) \cong \bigsqcup_{\la\in\cP_n^+}\B_n(\la)^{\oplus |{\tt L}^\la_{\mu\nu}|}.
\end{equation}
Hence we have $|{\tt L}^\la_{\mu\nu}|=f^\la_{\mu\nu}=f^\la_{\nu\mu}$ 
from \eqref{eq:P-function} and the linear independence of Schur $P$-polynomials $P_\la(x_1,\ldots,x_n)$'s.

Let us prove $f^\la_{\mu\nu}=|\tt{F}^\la_{\mu\nu}|$ by constructing a bijection 
\begin{equation}\label{eq:LRS recording}
\xymatrixcolsep{2pc}\xymatrixrowsep{2mm}\xymatrix@1{
{\tt L}^\la_{\mu\nu} \   \    \ar[r] \ar@{->}[r]  & \  \ {\tt F}^\la_{\mu\nu}  \\
 \  T  \ \ \  \ar[r] \ar@{|->}[r] &  \ \   Q_T.}
\end{equation}
Let $T \in {\tt L}^\la_{\mu\nu}$ be given.
Assume that $w_\rev(T)=u_1\cdots u_N$ where $N=|\nu|$. 
By Lemma \ref{lem:q(n)-lattice property}, there exists $\mu^{(m)}\in \cP_n^+$ for $1\leq m\leq N$  such that 
\begin{itemize}
\item[(i)] $(u_{N-m+1} \cdots u_N) \otimes L^\mu_n \equiv L^{\mu^{(m)}}_n$ and $\mu^{(N)}=\la$,

\item[(ii)] $\mu^{(m)}$ is obtained by adding a box in the $(n-u_m+1)$-st row of $\mu^{(m-1)}$.

\end{itemize}
Here we assume that $\mu^{(0)}=\mu$.
Recall that 
\begin{equation*}
w_\rev(T)=T^{(\ell(\nu))}\cdots T^{(1)},
\end{equation*}
where $T^{(k)}=T_{k,1}\cdots T_{k,\la_k}$ is a hook word for $1\leq k\leq \ell(\nu)$.
We define $Q_T$ to be a tableau of shifted shape $\la/\mu$ with entries in $\mc N$, where $\mu^{(m)}/\mu^{(m-1)}$ is filled with
\begin{equation}\label{eq:labeling of Q_T}
\begin{cases}
k', & \text{if $u_m$ belongs to $T^{(k)}\!\!\uparrow$},\\
k, & \text{if $u_m$ belongs to $T^{(k)}\!\!\downarrow$},\\
\end{cases}
\end{equation}
for some $1\leq k\leq \ell(\nu)$. 
In other words, the boxes in $Q_T$ corresponding to $T^{(k)}\!\!\uparrow$ are filled with $k'$ from right to left as a vertical strip 
and then those corresponding to $T^{(k)}\!\!\downarrow$ are filled with $k$ from left to right as a horizontal strip.

By construction, it is clear that $Q_T\in SST^+_{\mc N}(\la/\mu)$ with $c_{k'}(Q_T)+c_{k}(Q_T)=\nu_k$ for $1\leq k\leq \ell(\nu)$. Let $w(Q_T)=w_1\cdots w_N$. 
Since $T^{(k)}$ is a hook word for each $k$ and the rightmost letter, say $u_m$, in $T^{(k)}\!\!\downarrow$ is strictly smaller than the leftmost letter $u_{m+1}$ in $T^{(k)}\!\!\uparrow$, the entry $k$ in $Q_T$ corresponding to $u_m$ is located to the southeast of all $k'$'s in $Q_T$. So the conditions Definition \ref{Def:new LR tab}(1) and (2) are satisfied. 

It remains to check that $w(Q_T)$ satisfies the ``lattice property". 
Note that if we label $k$ and $k'$ in \eqref{eq:labeling of Q_T} 
as $k_j$ and $k'_j$, respectively when $u_m=T_{k,j}$, then it coincides with the labeling on the letters in $w(Q_T)$ given in Definition \ref{Def:w and w^*}(1). Now it is not difficult to see that the conditions Proposition \ref{prop:SSDT-checker}(1), (2), and (3) on $T$ implies the conditions Definition \ref{Def:new lattice property} (L1), (L2), and (L3), respectively. 
Therefore, $Q_T\in {\tt F}^\la_{\mu\nu}$.

Finally the correspondence $T \mapsto Q_T$ is injective and also reversible. Hence 
the map~\eqref{eq:LRS recording} is a bijection. This completes the proof. \qed

\begin{rem}
We see from Remark \ref{rem:stability} that $|{\tt L}^\la_{\mu\nu}|$ does not depend on $n$ for all sufficiently large $n$. Hence \eqref{eq:tensor decomp of mu nu} also implies the Schur $P$-positivity of the product $P_\mu P_\nu$.
\end{rem}

\begin{rem}
For $T\in {\tt L}^\la_{\mu\nu}$, let $\widehat{Q}_T$ be the tableau of shifted shape $\la/\mu$, which is defined in the same way as $Q_{T}$ in the proof of Theorem \ref{Thm:our-characterization} except that we fill $\mu^{(m)}/\mu^{(m-1)}$ with $m$ in \eqref{eq:labeling of Q_T} for $1\leq m\leq N$. Then the set $\{\,\widehat{Q}_T\,|\,T\in {\tt L}^\la_{\mu\nu}\,\}$ is equal to the one given in 
\cite[Theorem 4.13]{GJKKK14} to describe $f^\la_{\mu\nu}$. 
For example, 
\vskip 3mm
\begin{center}
\begin{tikzpicture}[baseline=0mm]
\def \hhh{4.4mm}
\def \vvv{5.4mm}
\draw[-] (\hhh*0,\vvv*0) rectangle (\hhh*1,\vvv*1);
\draw[-] (\hhh*1,\vvv*0) rectangle (\hhh*2,\vvv*1);
\draw[-] (\hhh*2,\vvv*0) rectangle (\hhh*3,\vvv*1);
\draw[-] (\hhh*1,-\vvv*1) rectangle (\hhh*2,\vvv*0);
\node at (\hhh*0.5,\vvv*0.5) {\small $3$};
\node at (\hhh*1.5,\vvv*0.5) {\small $3$};
\node at (\hhh*2.5,\vvv*0.5) {\small $4$};
\node at (\hhh*1.5,-\vvv*0.5) {\small $2$};
\node[left] at (-\hhh*0,\vvv*0) () {$T_1=$};
\node[right] at (\hhh*3,\vvv*0) {$\in {\tt L}^{(4,3,1)}_{(3,1)(3,1)}$};
\end{tikzpicture} 
\hskip 4mm
\begin{tikzpicture}[baseline=0mm]
\def \hhh{4.4mm}
\def \vvv{5.4mm}
\draw[-,black!20,fill=black!10] (\hhh*0,\vvv*0) rectangle (\hhh*1,\vvv*1);
\draw[-,black!20,fill=black!10] (\hhh*1,\vvv*0) rectangle (\hhh*2,\vvv*1);
\draw[-,black!20,fill=black!10] (\hhh*2,\vvv*0) rectangle (\hhh*3,\vvv*1);
\draw[-] (\hhh*3,\vvv*0) rectangle (\hhh*4,\vvv*1);
\draw[-,black!20,fill=black!10] (\hhh*1,-\vvv*1) rectangle (\hhh*2,\vvv*0);
\draw[-] (\hhh*2,-\vvv*1) rectangle (\hhh*3,\vvv*0);
\draw[-] (\hhh*3,-\vvv*1) rectangle (\hhh*4,\vvv*0);
\draw[-] (\hhh*2,-\vvv*2) rectangle (\hhh*3,-\vvv*1);
\node at (\hhh*3.5,\vvv*0.5) {\small $1$};
\node at (\hhh*2.5,-\vvv*0.5) {\small $2$};
\node at (\hhh*3.5,-\vvv*0.5)  {\small $3$};
\node at (\hhh*2.5,-\vvv*1.5) {\small $4$};
\node[left] at (-\hhh*0,-\vvv*0) () {$\widehat{Q}_{T_1}=$};
\end{tikzpicture} 
\hskip 12mm
\begin{tikzpicture}[baseline=0mm]
\def \hhh{4.4mm}
\def \vvv{5.4mm}
\draw[-,black!20,fill=black!10] (\hhh*0,\vvv*0) rectangle (\hhh*1,\vvv*1);
\draw[-,black!20,fill=black!10] (\hhh*1,\vvv*0) rectangle (\hhh*2,\vvv*1);
\draw[-,black!20,fill=black!10] (\hhh*2,\vvv*0) rectangle (\hhh*3,\vvv*1);
\draw[-] (\hhh*3,\vvv*0) rectangle (\hhh*4,\vvv*1);
\draw[-,black!20,fill=black!10] (\hhh*1,-\vvv*1) rectangle (\hhh*2,\vvv*0);
\draw[-] (\hhh*2,-\vvv*1) rectangle (\hhh*3,\vvv*0);
\draw[-] (\hhh*3,-\vvv*1) rectangle (\hhh*4,\vvv*0);
\draw[-] (\hhh*2,-\vvv*2) rectangle (\hhh*3,-\vvv*1);
\node at (\hhh*3.5,\vvv*0.5) {\small $1'$};
\node at (\hhh*2.5,-\vvv*0.5) {\small $1$};
\node at (\hhh*3.5,-\vvv*0.5)  {\small $1$};
\node at (\hhh*2.5,-\vvv*1.5) {\small $2$};
\node[left] at (-\hhh*0,-\vvv*0) {$Q_{T_1}=$};
\end{tikzpicture} 
\vskip 3mm
\begin{tikzpicture}[baseline=0mm]
\def \hhh{4.4mm}
\def \vvv{5.4mm}
\draw[-] (\hhh*0,\vvv*0) rectangle (\hhh*1,\vvv*1);
\draw[-] (\hhh*1,\vvv*0) rectangle (\hhh*2,\vvv*1);
\draw[-] (\hhh*2,\vvv*0) rectangle (\hhh*3,\vvv*1);
\draw[-] (\hhh*1,-\vvv*1) rectangle (\hhh*2,\vvv*0);
\node at (\hhh*0.5,\vvv*0.5) {\small $4$};
\node at (\hhh*1.5,\vvv*0.5) {\small $2$};
\node at (\hhh*2.5,\vvv*0.5) {\small $3$};
\node at (\hhh*1.5,-\vvv*0.5) {\small $3$};
\node[left] at (-\hhh*0,\vvv*0) () {$T_2=$};
\node[right] at (\hhh*3,\vvv*0) {$\in {\tt L}^{(4,3,1)}_{(3,1)(3,1)}$};
\end{tikzpicture} 
\hskip 4mm
\begin{tikzpicture}[baseline=0mm]
\def \hhh{4.4mm}
\def \vvv{5.4mm}
\draw[-,black!20,fill=black!10] (\hhh*0,\vvv*0) rectangle (\hhh*1,\vvv*1);
\draw[-,black!20,fill=black!10] (\hhh*1,\vvv*0) rectangle (\hhh*2,\vvv*1);
\draw[-,black!20,fill=black!10] (\hhh*2,\vvv*0) rectangle (\hhh*3,\vvv*1);
\draw[-] (\hhh*3,\vvv*0) rectangle (\hhh*4,\vvv*1);
\draw[-,black!20,fill=black!10] (\hhh*1,-\vvv*1) rectangle (\hhh*2,\vvv*0);
\draw[-] (\hhh*2,-\vvv*1) rectangle (\hhh*3,\vvv*0);
\draw[-] (\hhh*3,-\vvv*1) rectangle (\hhh*4,\vvv*0);
\draw[-] (\hhh*2,-\vvv*2) rectangle (\hhh*3,-\vvv*1);
\node at (\hhh*3.5,\vvv*0.5) {\small $3$};
\node at (\hhh*2.5,-\vvv*0.5) {\small $1$};
\node at (\hhh*3.5,-\vvv*0.5)  {\small $4$};
\node at (\hhh*2.5,-\vvv*1.5) {\small $2$};
\node[left] at (-\hhh*0,-\vvv*0) () {$\widehat{Q}_{T_2}=$};
\end{tikzpicture} 
\hskip 12mm
\begin{tikzpicture}[baseline=0mm]
\def \hhh{4.4mm}
\def \vvv{5.4mm}
\draw[-,black!20,fill=black!10] (\hhh*0,\vvv*0) rectangle (\hhh*1,\vvv*1);
\draw[-,black!20,fill=black!10] (\hhh*1,\vvv*0) rectangle (\hhh*2,\vvv*1);
\draw[-,black!20,fill=black!10] (\hhh*2,\vvv*0) rectangle (\hhh*3,\vvv*1);
\draw[-] (\hhh*3,\vvv*0) rectangle (\hhh*4,\vvv*1);
\draw[-,black!20,fill=black!10] (\hhh*1,-\vvv*1) rectangle (\hhh*2,\vvv*0);
\draw[-] (\hhh*2,-\vvv*1) rectangle (\hhh*3,\vvv*0);
\draw[-] (\hhh*3,-\vvv*1) rectangle (\hhh*4,\vvv*0);
\draw[-] (\hhh*2,-\vvv*2) rectangle (\hhh*3,-\vvv*1);
\node at (\hhh*3.5,\vvv*0.5) {\small $1$};
\node at (\hhh*2.5,-\vvv*0.5) {\small $1'$};
\node at (\hhh*3.5,-\vvv*0.5)  {\small $2$};
\node at (\hhh*2.5,-\vvv*1.5) {\small $1$};
\node[left] at (-\hhh*0,-\vvv*0) {$Q_{T_2}=$};
\end{tikzpicture} 
\end{center}
\end{rem}

\subsection{Stembridge's description of $f^\la_{\mu\nu}$}

\begin{df}\label{df:St lattice property}
Let $w=w_1 \cdots w_N$ be a word in $\W_{\mc N}$
and $w_\rev$ be the reverse word of $w$.
Let $\widehat{w}$ be the word obtained from $w$ by replacing $k$ by $(k+1)'$ and $k'$ by $k$ for each $k \geq 1$.
Suppose that $w \widehat{w}_\rev = a_1\cdots a_{2N}$, and let ${m}_k(i)=c_k(a_{1}  \cdots  a_i)$ for $k \geq 1$ and $0 \leq i \leq 2N$. 
Then we say that $w$ satisfies the {\em lattice property} 
if 
\begin{equation}\label{Eqn:the lattice property}
\text{$m_{k+1}(i)=m_{k}(i)$ implies $|a_{i+1}|\neq k+1$ for $k\geq 1$ and $i\geq 0$}.
\end{equation}
Here we assume that $m_k(0)=0$.
\end{df}

\begin{df}
\label{Def:Stembridge's description}
For $\la, \mu, \nu \in \cP^+$, let 
${\tt LRS}^\la_{\mu\nu}$ be the set of tableaux $Q$ such that 
\begin{enumerate}
\item[(1)] $Q\in SST^+_{\mc N}(\la/\mu)$ with $c_{k}(Q)+c_{k'}(Q)=\nu_k$ for $k\geq 1$,

\item[(2)] for $k\geq 1$, if $x$ is the rightmost letter in $w(Q)$ with $|x|=k$, then $x=k$,

\item[(3)] $w(Q)$ satisfies the lattice property in Definition \ref{df:St lattice property}.
\end{enumerate}
We call ${\tt LRS}^\la_{\mu\nu}$ the set of {\em Littlewood-Richardson-Stembridge tableaux}.
\end{df}

\begin{thm}\label{Stembridge's result}{\rm (\cite[Theorem 8.3]{Ste})}
For $\la, \mu, \nu \in \cP^+$, we have 
$$
f^\la_{\mu\nu} = \left|{\tt LRS}^\la_{\mu\nu}\right|,
$$
that is, the shifted LR coefficient $f^\la_{\mu\nu}$ is equal to the number of tableaux in ${\tt LRS}^\la_{\mu\nu}$.
\end{thm}

\begin{thm}\label{thm:F=LRS}
For $\lambda, \mu, \nu \in \cP^+$, we have
\begin{displaymath}
\tt{F}^\la_{\mu\nu}={\tt LRS}^\la_{\mu\nu}.
\end{displaymath}
\end{thm}
\pf Since  
Definition~\ref{Def:new LR tab}(1) and (2) are the same as 
Definition~\ref{Def:Stembridge's description}(1) and (2), respectively,
it suffices to show that 
for any $Q \in SST^+_{\mc N}(\lambda/\mu)$, 
$w:=w(Q)$ satisfies the ``lattice property'' in Definition~\ref{Def:new lattice property}  
if and only if $w$ satisfies the lattice property in Definition~\ref{df:St lattice property}.
We assume that $N = |\nu|$,
$w=w_1 \cdots w_N$,
$w^* = w^*_1  \cdots w^*_N$,
and
$w \widehat{w}_{\rev} =  a_1 \cdots a_{2N}$.
\vskip 1mm

Suppose that 
$w$ satisfies the ``lattice property'' in Definition~\ref{Def:new lattice property}.  
We use induction on $1 \leq i \leq 2N$ to show that $a_1 \cdots a_{2N}$ satisfies~\eqref{Eqn:the lattice property}.
We first observe from (L1) that $a_1 = $  $1$ or $1'$, and  $a_1$ satisfies \eqref{Eqn:the lattice property} since $m_{k}(0)=0$ for all $k\geq 1$.

We now assume that  
$a_1 \cdots a_i$ for some $1 \leq i < 2N$ satisfies \eqref{Eqn:the lattice property}.
Suppose for the sake of contradiction that 
$m_{k+1}(i) = m_{k}(i) =m$ and $|a_{i+1}| = k+1$ for some $k \geq 1$. 
Here $m>0$ by (L1).  
By induction hypothesis, there exist $s < t \leq i$  such that 
$a_s = k$ with $m_k(s)=m$  and  $a_t = k+1$ with $m_{k+1}(t)=m$.
Note that for each $k \geq 1$
\begin{equation}\label{Eq:strict w}
c_k(w)+c_{k'}(w) > c_{k+1}(w) + c_{(k+1)'}(w),
\end{equation}
which implies that the number of $k$'s in $w \widehat{w}_{\rev}$ is greater than the number  of $k+1$'s in $w \widehat{w}_{\rev}$.
So we can choose an integer $u>i+1$ such that $a_u=k$ and $m_k(u)=m+1$.
We now consider the following four cases:
\vskip 1mm
{\em Case 1}. Let $1 \leq s < t < i+1 \leq N$.
In this case $w^*_s = k_m$, $w^*_t=k+1_m$ and $w^*_{i+1}=k+1_{M}$ for some $M\geq m+1$.
(i) If $u \leq N$, then  
we have $(w^*_t, w^*_{i+1},w^*_u) = (k+1_m, k+1_{M}, k_{m+1})$, which contradicts (L2).
%
(ii) If $N < u < 2N-s+1$, then $a_{2N-u+1}=w_{2N-u+1}=k'$ ($s<2N-u+1\leq N$) but  no $k$ occurs in $w_{s+1}\cdots w_N$  
which contradicts Definition~\ref{Def:new LR tab}(2).
%
(iii) If $2N-s+1 < u \leq 2N$, then 
we have $(w^*_{2N-u+1}, w^*_s,w^*_t) = (k_{m+1}, k_{m}, k+1_{m})$, which contradicts (L3).
\vskip 1mm
{\em Case 2}. Let $1 \leq s < t \leq N < i+1 \leq 2N$.
In this case $w^*_s = k_m$ and $w^*_t=k+1_m$. Since $w_{2N-u+1}=k'$ ($2N-u+1 < N$), we have $s\neq 2N-u+1$.
(i) If $s < 2N-u+1$, then
we have $w_{2N-u+1}=k'$ but no $k$ in $w_{s+1}\cdots w_n$ since $m_{k}(i)=m$, which contradicts Definition~\ref{Def:new LR tab}(2).
(ii) If $2N-u+1 < s$, then
we have $(w^*_{2N-u+1}, w^*_{s},w^*_{t}) = (k_{m+1}, k_{m}, k+1_{m})$, which contradicts (L3).
\vskip 1mm
{\em Case 3}. Let $1 \leq s \leq N < t < i+1 \leq 2N$.
In this case $w^*_s = k_m$, $w^*_{2N-t+1}=k+1_m$.
If $a_{i+1}=(k+1)'$, then $a_{2N-i}=k$ but it is impossible from the assumption $m_k(i) = m$.
So $a_{i+1}= k+1$ and $w^*_{2N-i}=k+1_{m+1}$.
(i) If $s < 2N-t+1$, then we have $w_{2N-t+1}=(k+1)'$ ($2N-t+1\leq N$) but no $k+1$ in $w_{2N-t+2}\cdots w_N$ since $m_{k+1}(i)=m$, which contradicts Definition~\ref{Def:new LR tab}(2).
(ii) If $2N-t+1 < s$, then by \eqref{Eq:strict w}
there is an integer $v > u$ such that $a_v = k$ and $m_k(v)=m+2$.
So we have $(w^*_{2N-v+1}, w^*_{2N-u+1},w^*_{2N-i}) = (k_{m+2}, k_{m+1}, k+1_{m+1})$, which contradicts (L3).
\vskip 1mm
{\em Case 4}. Let $N < s < t < i+1 \leq 2N$.
In this case $w^*_{2N-s+1} = k_m$ and $w^*_{2N-t+1}=k+1_m$.
(i) If $a_{i+1}=k+1$, then $w_{2N-i}=(k+1)'$ and $w^*_{2N-i}=k+1_{m+1}$.   
By \eqref{Eq:strict w}
there is an integer $v > u$ such that $a_v = k$ and $m_k(v)=m+2$.
So we have $(w^*_{2N-v+1}, w^*_{2N-u+1},w^*_{2N-i}) = (k_{m+2}, k_{m+1}, k+1_{m+1})$, which contradicts (L3). 
(ii) If $a_{i+1}=(k+1)'$, then $w_{2N-i}=k$ and $w^*_{2N-i}=k_{M}$ for some $M<m$. So we have $(w^*_{2N-u+1}, w^*_{2N-i},w^*_{2N-t+1}) = (k_{m+1}, k_{M}, k+1_{m})$, which contradicts (L3).
\vskip 1mm

Conversely, we assume that $w$ satisfies the lattice property in Definition~\ref{df:St lattice property}. 
We first claim that $w$ satisfies (L1). 
Given $k\geq 1$, let $w^*_i = k_1$ for some $1 \leq i \leq N$. 
If $w^*_j = k+1_1$ for some $1 \leq j < i$, then it follows that 
$m_{k}(j-1)=m_{k+1}(j-1)=0$ and $a_j =  k+1$, which contradicts~\eqref{Eqn:the lattice property}. Hence $w$ satisfies (L1).

Next, we claim that $w$ satisfies (L2).
Suppose  that   
there is a triple $(w^*_s, w^*_u, w^*_t) = (k+1_i, k+1_{j}, k_{i+1})$
for some $k\geq 1$, $i<j$, and $1 \leq s < u <  t \leq N$. We may assume that $j=i+1$.
Since $w^*_s = k+1_i$ is placed to the left of $w^*_u = k+1_{i+1}$, 
it follows from Definition~\ref{Def:w and w^*} that $a_s= k+1$, and
from Definition~\ref{Def:Stembridge's description}(2) that $a_u = k+1$.
Since $w$ satisfies the lattice property, 
there is a positive integer $v < s$ such that $w^*_v=k_{i}$,
i.e., $a_v=k$ and $m_{k}(v)=i$ for some $v<s$.
We have $m_{k+1}(u-1) = m_{k}(u-1)=i$ and $a_u = k+1$, a contradiction. So $w$ satisfies (L2).

Finally, we claim that $w$ satisfies (L3).
Suppose for the sake of contradiction that $(w^*_s, w^*_u, w^*_t)=(k_{j+1}, k_i, k+1_j)$ for some $k\geq 1$, $i\leq j$, and $1 \leq s<u<t \leq N$. We may assume that $i=j$.
Since $w^*_s = k_{j+1}$ is placed to the left of $w^*_u = k_j$, it follows that $a_s = k'$.
We consider four cases depending on the primedness of $a_u$ and $a_t$ as follows:

{\em Case 1}. Let $a_u=k'$ and $a_t=(k+1)'$.
It follows that $a_{2N-u+1}=k$ ($m_{k}(2N-u+1)=j$) and $a_{2N-t+1}=k+1$ ($m_{k}(2N-u+1)=j)$.
So we have $m_{k+1}(2N-t)=m_{k}(2N-t)=j-1$ and $a_{2N-t+1} =k+1$, as desired.

{\em Case 2}. Let $a_u=k'$ and $a_t=k+1$.
It follows that $a_{2N-u+1}=k$, $m_{k}(2N-u+1)=j$ and $m_{k+1}(t)=j$.
Since $t < 2N-u+1$, we have $m_{k}(t) < m_{k+1}(t)$.
So there is an integer $0 \leq \hat{t} < t$ such that 
$m_{k}(\hat{t}) = m_{k+1}(\hat{t}) <j$ and $a_{\hat{t}+1}=k+1$, as desired.

{\em Case 3}. Let $a_u = k$ and $a_t = (k+1)'$. 
It follows that $a_{2N-t+1}=k+1$ and $m_{k+1}(2N-t+1)=j$.
From $m_{k}(2N-s+1)=j+1$ and $2N-t+1<2N-s+1$ we have $m_{k}(2N-t+1) = m_{k+1}(2N-t+1)=j$.
If there is another $k+1$ between $w^*_t$ and $w^*_u$,
then we obtain the desired contradiction.
Otherwise, 
$m_{k}(2N-u) = m_{k+1}(2N-u)$ and $a_{2N-u+1} = (k+1)'$, as desired.

{\em Case 4}.
Let $a_u = k$ and $a_t = k+1$. 
From $a_s=k'$ ($w^*_s=k_{j+1}$) it follows that $m_{k}(2N-u) = j$.
If $m_{k+1}(2N-u)=j$, from $a_{2N-u+1} = (k+1)'$ we get a contradiction.
If $m_{k+1}(2N-u)>j$, by choosing the smallest integer $\hat{t}>t$ such that $m_{k+1}(\hat{t})=j+1$ this leads to a contradiction.
\qed
\vskip 3mm

Indeed, we have shown in the proof of Theorem \ref{thm:F=LRS} that
\begin{cor}\label{cor:equivalence of lattice conditions}
Let $w\in \W_{\mc N}$ be such that 
\begin{itemize}
\item[(1)]  $(c_{k}(Q)+c_{k'}(Q))_{k\geq 1}\in \cP^+$,

\item[(2)] for $k\geq 1$, if $x$ is the rightmost letter in $w$ with $|x|=k$, then $x=k$.
\end{itemize}
Then $w$ satisfies the ``lattice property" in Definition \ref{Def:new lattice property} if and only if $w$ satisfies the lattice property in Definition \ref{df:St lattice property}.
\end{cor}

\begin{rem}
A bijection from ${\tt LRS}^{\la}_{\mu \nu}$ to ${\tt L}^\la_{\mu \nu}$ is also given in \cite[Theorem 4.7]{CNO14}, which coincides with the inverse of the map $T \mapsto Q_T$ in \eqref{eq:LRS recording} (see also the remarks in \cite[p.82]{CNO14}).
The proof of \cite[Theorem 4.7]{CNO14} use insertion schemes for two versions of  semistandard decomposition tableaux and another combinatorial model for $f^{\la}_{\mu\nu}$ by Cho \cite{Cho} as an intermediate object between ${\tt LRS}^{\la}_{\mu \nu}$ and ${\tt L}^\la_{\mu \nu}$.

On the other hand, we prove more directly that the map $T \mapsto Q_T$ in \eqref{eq:LRS recording} is a bijection from ${\tt L}^\la_{\mu \nu}$ to $\tt{LRS}^\la_{\mu\nu}$ by using a new characterization of the lattice property in Theorem \ref{thm:F=LRS}.
\end{rem}

\section{Schur $P$-expansions of skew Schur functions}

\subsection{The Schur $P$-expansion of $s_{\la/\delta_r}$}

For $r\geq 0$, let us denote by $\delta_r$ the partition $(r,r-1,\ldots,1)$ if $r \geq 1$,
and $(0)$ if $r=0$. We fix a non-negative integer $r$.

Let $\la\in \cP$ be such that $D_{\delta_{r}}\subseteq D_\la \subseteq D_{((r+1)^{r+1})}$.
Here $((r+1)^{r+1})$ means the rectangular partition $(r+1,\ldots,r+1)$ with length $r+1$.
For instance, the diagram 
\begin{displaymath}
\begin{tikzpicture}[baseline=0mm]
\def \hhh{3.8mm}    
\def \vvv{4.5mm}    
\node[left] at (-\hhh*0,-\vvv*1) {$D_{(5,4,4,4,2)/{\delta_4}}=$};
\draw[-,black!10] (\hhh*4,-\vvv*1) rectangle (\hhh*5,-\vvv*0);
\draw[-,black!10] (\hhh*4,-\vvv*2) rectangle (\hhh*5,-\vvv*1);
\draw[-,black!10] (\hhh*4,-\vvv*3) rectangle (\hhh*5,-\vvv*2);
\draw[-,black!10] (\hhh*2,-\vvv*4) rectangle (\hhh*3,-\vvv*3);
\draw[-,black!10] (\hhh*3,-\vvv*4) rectangle (\hhh*4,-\vvv*3);
\draw[-,black!10] (\hhh*4,-\vvv*4) rectangle (\hhh*5,-\vvv*3);
\draw[-,black!20,fill=black!10] (\hhh*0,0) rectangle (\hhh*1,\vvv*1);
\draw[-,black!20,fill=black!10] (\hhh*1,0) rectangle (\hhh*2,\vvv*1);
\draw[-,black!20,fill=black!10] (\hhh*2,0) rectangle (\hhh*3,\vvv*1);
\draw[-,black!20,fill=black!10] (\hhh*3,0) rectangle (\hhh*4,\vvv*1);
\draw[-] (\hhh*4,0) rectangle (\hhh*5,\vvv*1);
\draw[-,black!20,fill=black!10] (\hhh*0,-\vvv*1) rectangle (\hhh*1,\vvv*0);
\draw[-,black!20,fill=black!10] (\hhh*1,-\vvv*1) rectangle (\hhh*2,\vvv*0);
\draw[-,black!20,fill=black!10] (\hhh*2,-\vvv*1) rectangle (\hhh*3,\vvv*0);
\draw[-] (\hhh*3,-\vvv*1) rectangle (\hhh*4,\vvv*0);
\draw[-,black!20,fill=black!10] (\hhh*0,-\vvv*2) rectangle (\hhh*1,-\vvv*1);
\draw[-,black!20,fill=black!10] (\hhh*1,-\vvv*2) rectangle (\hhh*2,-\vvv*1);
\draw[-] (\hhh*2,-\vvv*2) rectangle (\hhh*3,-\vvv*1);
\draw[-] (\hhh*3,-\vvv*2) rectangle (\hhh*4,-\vvv*1);
\draw[-,black!20,fill=black!10] (\hhh*0,-\vvv*3) rectangle (\hhh*1,-\vvv*2);
\draw[-] (\hhh*1,-\vvv*3) rectangle (\hhh*2,-\vvv*2);
\draw[-] (\hhh*2,-\vvv*3) rectangle (\hhh*3,-\vvv*2);
\draw[-] (\hhh*3,-\vvv*3) rectangle (\hhh*4,-\vvv*2);
\draw[-] (\hhh*0,-\vvv*4) rectangle (\hhh*1,-\vvv*3);
\draw[-] (\hhh*1,-\vvv*4) rectangle (\hhh*2,-\vvv*3);
\end{tikzpicture}
\end{displaymath} 
is contained in $D_{(5^5)}$.

It is shown in~\cite{AS,Dew}
that the skew Schur function $s_{\lambda/\delta_k}$ has a non-negative integral expansion in terms of Schur $P$-functions
\begin{equation}\label{Eq:s-to-P}
s_{\lambda/\delta_r} = \sum_{\nu\in\cP^+} a_{\lambda/\delta_r\,\nu} \hskip 0.5mm P_\nu, 
\end{equation}
together with a combinatorial description of $a_{\lambda/\delta_r\,\nu}$.
Moreover it is shown that these skew Schur functions are the only ones (up to rotation of shape by $180^\circ$), which have Schur $P$-positivity. In this section, we give a new simple description of $a_{\lambda/\delta_r\,\nu}$ using $\q(n)$-crystals.

First we consider a $\mf q(n)$-crystal structure on $B_n(\la/\delta_r)$, which is a slight generalization of \cite[Example 2.10(d)]{GJKKK15}.

\begin{prop}\label{Prop:B-q(n)-crystal}
Let $\lambda \in \cP_n$ be such that $D_{\delta_{r}} \subseteq D_\lambda \subseteq D_{(r+1)^{r+1}}$.
Then the $\gl(n)$-crystal $B_n(\la/\delta_r)$ as a subset of $\W_{[n]}$ together with ${\bf 0}$ is invariant under $\te_{\ov{1}}$ and $\tf_{\ov{1}}$. Hence $B_n(\la/\delta_r)$ is a $\q(n)$-crystal.
\end{prop}
\pf
Let $N = |\la|-|\delta_r|$.
For $T \in B_n(\la/\delta_r)$, let $w(T) = w_1 \cdots w_N$. 
Recall that $T$ is identified with $w(T)$ in $\W_{[n]}$.
Here we call the box in $D_{\lambda/\delta_r}$ containing $w_i$ the $w_i$-box, and call
the set of boxes $(x,r-x+2) \in D_{\lambda/\delta_r}$ for $1 \leq x \leq r+1$
the main anti-diagonal of $D_{\lambda/\delta_r}$.
 
Suppose that $\tf_{\overline{1}} w(T) \neq 0$.
There exists $1\leq i \leq N-1$ such that 
$w_i = 1$ and $w_j \neq 1, 2$ for all $i < j \leq N$, and   
\begin{displaymath}
\tf_{\ov 1} (w_1 \cdots w_{i-1} \,  1 \, w_{i+1} \cdots w_{N}) = w_1 \cdots w_{i-1} \ 2  \ w_{i+1} \cdots w_{N},
\end{displaymath}
by the tensor product rule \eqref{eq:tensor product of qn-crystals}.
We first observe that the entry $1$ in $T$  
can be placed only on the main anti-diagonal in $D_{\lambda/\delta_r}$. 
If there is a box in $D_{\lambda/\delta_r}$ below the $w_i$-box, 
then it corresponds to $w_j$ for some $j > i$, and hence its entry is greater than 2.
Moreover, if there is a box in $D_{\lambda/\delta_r}$ to the right of the $w_i$-box, then its entry is greater than $1$ since it is not on the main anti-diagonal.
So we conclude that there exists $T' \in SST_{[n]}(\la/\delta_r)$
such that $w(T') = \tf_{\ov 1} w(T)$.
\vskip 1mm
Suppose that $\te_{\overline{1}} w(T) \neq 0$.
There exists $1\leq i \leq N-1$ such that 
$w_i = 2$ and $w_j \neq 1, 2$ for all $i < j \leq N$, and
\begin{equation}\label{eq:aux-1}
\te_{\ov 1}(w_1 \cdots w_{i-1} \,  2 \,  w_{i+1} \cdots w_{N}) = w_1 \cdots w_{i-1} \ 1  \ w_{i+1} \cdots w_{N},
\end{equation}
by the tensor product rule \eqref{eq:tensor product of qn-crystals}.
If the $w_i$-box is not on the main anti-diagonal, then
the $w_{i+1}$-box is placed to the left of the $w_i$-box.
Then the $w_{i+1}$-box is filled with $1$ or $2$, which contradicts \eqref{eq:aux-1}.
So the $w_i$-box is on the main anti-diagonal,
and thus $\tf_{\ov 1} w(T)=w(T')$ for some $T'\in B_n(\lambda/\delta_r)$. This completes the proof.
\qed
\vskip 3mm

\begin{cor}
Under the above hypothesis, the skew Schur function $s_{\la/\delta_r}$ is Schur $P$-positive.
\end{cor}
\pf Since $B_n(\la/\delta_r)$ is a $\mf q(n)$-crystal, the skew Schur polynomial $s_{\la/\delta_r}(x_1,\ldots,x_n)$ is a non-negative integral linear combination of $P_\nu(x_1,\ldots,x_n)$. Then we apply  Remark \ref{rem:stability}. 
\qed

\begin{df}\label{Def:Set-A}
Let $\la\in \cP$ be such that $D_{\delta_{r}} \subseteq D_\la \subseteq D_{((r+1)^{r+1})}$ and $\nu\in \cP^+$.
Let 
${\tt A}_{\la/\delta_r\,\nu}$ be the set of tableaux $Q$ such that 
\begin{itemize}
\item[(1)] $Q\in SST_{[r+1]}^+(\nu)$ with $c_{k}(Q)=\la_{r-k+2}-k+1$ for $1\leq k\leq r+1$,

\item[(2)] for $1\leq k\leq r$ and $1\leq i\leq N$,  
\begin{equation*}\label{Rule:S-to-P}
m_{k}(i)\leq m_{k+1}(i) +1,
\end{equation*} 
where $w_\rev(Q)=w_1\cdots w_N$ and 
$m_k(i)=c_k(w_1\cdots w_i)$.
\end{itemize}
\end{df}

Then we have the following combinatorial description of $a_{\la/\delta_r\,\nu}$. 
\begin{thm}\label{Thm:Main-sect5}
For $\la\in \cP$ with $D_{\delta_{r}} \subseteq D_\la \subseteq D_{((r+1)^{r+1})}$ 
and $\nu\in \cP^+$, we have
\begin{equation*}
a_{\la/\delta_r\,\nu}=\left|{\tt A}_{\la/\delta_r\,\nu}\right|.
\end{equation*}
\end{thm}
\pf  Choose $n$ such that $\la, \nu\in \cP_n^+$. 
We may assume that $\lambda_1 = \ell(\la) =  r+1$.
Let 
\begin{equation}\label{eq:L_{la del}}
{\tt L}_{\la/\delta_r\,\nu}=\{\,T \in B_n(\lambda/\delta_r) \,|\, T\equiv L^\nu\,\}.
\end{equation} 
By Proposition \ref{Prop:B-q(n)-crystal}, we have
\begin{equation}\label{eq:decomp of stair}
B_n(\la/\delta_r) \cong \bigsqcup_{\nu \in \cP_n^+} \B_n(\nu)^{\oplus |{\tt L}_{\la/\delta_r\,\nu}|}.
\end{equation}
By linear independence of $P_{\nu}(x_1,\ldots,x_n)$'s for $\nu\in\cP_n^+$, we have $a_{\la/\delta_r\,\nu}=|{\tt L}_{\la/\delta_r\,\nu}|.$ 

Let us construct a bijection
\begin{equation}\label{eq:Q_T for a}
\xymatrixcolsep{2pc}\xymatrixrowsep{0pc}\xymatrix{
\makebox[3em][c]{${\tt L}_{\la/\delta_r \, \nu}$} \ \ar[r] \ar@{->}[r]  & \  \makebox[3em][c]{${\tt A}_{\la/\delta_r \, \nu}$}  \\
\makebox[3em][c]{$T$}    \ \ar[r] \ar@{|->}[r] &    \makebox[3em][c]{$Q_T$}}
\end{equation}
as follows. 
Let $T \in {\tt L}_{\la/\delta_r \, \nu}$ be given.
Suppose that $w(T)=u_1\cdots u_N$, where $N=|\nu|$. 
By Lemma \ref{lem:q(n)-lattice property}, 
there exists $\nu^{(m)}\in \cP_n^+$ for $1\leq m\leq N$ such that $u_{N-m+1} \cdots u_N \equiv L^{\nu^{(m)}}$, where $\nu^{(1)}=(1)$, $\nu^{(N)} = \nu$, and
$\nu^{(m)}$ is obtained by adding a box in the $(n-u_m+1)$-st row of $\nu^{(m-1)}$ for $1\leq m\leq N$ with $\nu^{(0)}=\emptyset$. 

Note that $w_\rev(T)=T^{(r+1)}\cdots T^{(1)}$, 
where $T^{(l)}=T_{l,1}\cdots T_{l,\la_l-r-1+l}$ is a weakly increasing word corresponding to the $l$-th row of $T$ for $1\leq l\leq r+1$.
Let $Q_T$ be a tableau of shifted shape $\nu$ with entries in $\mathbb{N}$, where $\nu^{(m)}/\nu^{(m-1)}$ is filled with $r+2-l$
if $u_m$ occurs in $T^{(l)}$,
for some $1\leq l\leq r+1$. 
Note that the boxes in $Q_T$ corresponding to $T^{(l)}$ are filled with $r+2-l$ as a horizontal strip. So $Q_T$ satisfies the condition Definition \ref{Def:Set-A}(1).

For each $k\geq 1$, let us enumerate the letter $k$'s in $Q_T$ from southwest to northeast like $k_1,k_2,\ldots$.
Since $T\in SST_n(\la/\delta_r)$, we see that the entry $k_i$ in $Q_T$ corresponds to $T_{l,i}$  for $i\geq 1$, where $l=r+2-k$, and moreover $(k+1)_i$ is located in the southwest of $k_{i+1}$ for $i\geq 2$. This implies the condition Definition \ref{Def:Set-A}(2), and hence $Q_T\in {\tt A}_{\la/\delta_r \, \nu}$.

Finally, one can check that correspondence $T \mapsto Q_T$ is a bijection. 
\qed

\begin{exm}
Let $\lambda = (5,5,4,3,1)$ with $D_\la \subseteq D_{(5^5)}$ and $n=7$. 
For $\nu = (4,3,1)$,  
we have ${\tt L}_{\la/\delta_4\, \nu}=\{\,T_1, T_2\,\}$
and ${\tt A}_{\la/\delta_4\,\nu}=\{\,Q_{T_1}, Q_{T_2}\,\}$ as follows.
\begin{displaymath}
\begin{tikzpicture}[baseline=0mm]
\def \hhh{4mm}
\def \vvv{4.8mm}
\draw[-,black!20,fill=black!10] (\hhh*0,\vvv*0) rectangle (\hhh*1,\vvv*1);
\draw[-,black!20,fill=black!10] (\hhh*1,\vvv*0) rectangle (\hhh*2,\vvv*1);
\draw[-,black!20,fill=black!10] (\hhh*2,\vvv*0) rectangle (\hhh*3,\vvv*1);
\draw[-,black!20,fill=black!10] (\hhh*3,\vvv*0) rectangle (\hhh*4,\vvv*1);
\draw[-] (\hhh*4,\vvv*0) rectangle (\hhh*5,\vvv*1);
\draw[-,black!20,fill=black!10] (\hhh*0,-\vvv*1) rectangle (\hhh*1,\vvv*0);
\draw[-,black!20,fill=black!10] (\hhh*1,-\vvv*1) rectangle (\hhh*2,\vvv*0);
\draw[-,black!20,fill=black!10] (\hhh*2,-\vvv*1) rectangle (\hhh*3,\vvv*0);
\draw[-] (\hhh*3,-\vvv*1) rectangle (\hhh*4,\vvv*0);
\draw[-] (\hhh*4,-\vvv*1) rectangle (\hhh*5,\vvv*0);
\draw[-,black!20,fill=black!10] (\hhh*0,-\vvv*2) rectangle (\hhh*1,-\vvv*1);
\draw[-,black!20,fill=black!10] (\hhh*1,-\vvv*2) rectangle (\hhh*2,-\vvv*1);
\draw[-] (\hhh*2,-\vvv*2) rectangle (\hhh*3,-\vvv*1);
\draw[-] (\hhh*3,-\vvv*2) rectangle (\hhh*4,-\vvv*1);
\draw[-,black!20,fill=black!10] (\hhh*0,-\vvv*3) rectangle (\hhh*1,-\vvv*2);
\draw[-] (\hhh*1,-\vvv*3) rectangle (\hhh*2,-\vvv*2);
\draw[-] (\hhh*2,-\vvv*3) rectangle (\hhh*3,-\vvv*2);
\draw[-] (\hhh*0,-\vvv*4) rectangle (\hhh*1,-\vvv*3);
\node at (\hhh*4.5,\vvv*0.5) {\small $6$};
\node at (\hhh*3.5,-\vvv*0.5) {\small $5$};
\node at (\hhh*4.5,-\vvv*0.5) {\small $7$};
\node at (\hhh*2.5,-\vvv*1.5) {\small $6$};
\node at (\hhh*3.5,-\vvv*1.5) {\small $6$};
\node at (\hhh*1.5,-\vvv*2.5) {\small $7$};
\node at (\hhh*2.5,-\vvv*2.5) {\small $7$};
\node at (\hhh*0.5,-\vvv*3.5) {\small $7$};
%
%
\node[left] at (-\hhh*0.2,-\vvv*1) {$T_1=$};
\end{tikzpicture} 
\hskip 3mm
\begin{tikzpicture}[baseline=0mm]
\def \hhh{4mm}
\def \vvv{4.8mm}
\draw[-,black!20,fill=black!10] (\hhh*0,\vvv*0) rectangle (\hhh*1,\vvv*1);
\draw[-,black!20,fill=black!10] (\hhh*1,\vvv*0) rectangle (\hhh*2,\vvv*1);
\draw[-,black!20,fill=black!10] (\hhh*2,\vvv*0) rectangle (\hhh*3,\vvv*1);
\draw[-,black!20,fill=black!10] (\hhh*3,\vvv*0) rectangle (\hhh*4,\vvv*1);
\draw[-] (\hhh*4,\vvv*0) rectangle (\hhh*5,\vvv*1);
\draw[-,black!20,fill=black!10] (\hhh*0,-\vvv*1) rectangle (\hhh*1,\vvv*0);
\draw[-,black!20,fill=black!10] (\hhh*1,-\vvv*1) rectangle (\hhh*2,\vvv*0);
\draw[-,black!20,fill=black!10] (\hhh*2,-\vvv*1) rectangle (\hhh*3,\vvv*0);
\draw[-] (\hhh*3,-\vvv*1) rectangle (\hhh*4,\vvv*0);
\draw[-] (\hhh*4,-\vvv*1) rectangle (\hhh*5,\vvv*0);
\draw[-,black!20,fill=black!10] (\hhh*0,-\vvv*2) rectangle (\hhh*1,-\vvv*1);
\draw[-,black!20,fill=black!10] (\hhh*1,-\vvv*2) rectangle (\hhh*2,-\vvv*1);
\draw[-] (\hhh*2,-\vvv*2) rectangle (\hhh*3,-\vvv*1);
\draw[-] (\hhh*3,-\vvv*2) rectangle (\hhh*4,-\vvv*1);
\draw[-,black!20,fill=black!10] (\hhh*0,-\vvv*3) rectangle (\hhh*1,-\vvv*2);
\draw[-] (\hhh*1,-\vvv*3) rectangle (\hhh*2,-\vvv*2);
\draw[-] (\hhh*2,-\vvv*3) rectangle (\hhh*3,-\vvv*2);
\draw[-] (\hhh*0,-\vvv*4) rectangle (\hhh*1,-\vvv*3);
\node at (\hhh*4.5,\vvv*0.5) {\small $5$};
\node at (\hhh*3.5,-\vvv*0.5) {\small $6$};
\node at (\hhh*4.5,-\vvv*0.5) {\small $6$};
\node at (\hhh*2.5,-\vvv*1.5) {\small $6$};
\node at (\hhh*3.5,-\vvv*1.5) {\small $7$};
\node at (\hhh*1.5,-\vvv*2.5) {\small $7$};
\node at (\hhh*2.5,-\vvv*2.5) {\small $7$};
\node at (\hhh*0.5,-\vvv*3.5) {\small $7$};
%
%
\node[left] at (-\hhh*0.2,-\vvv*1) {$T_2=$};
\end{tikzpicture} 
\hskip 8mm
\begin{tikzpicture}[baseline=4mm]
\def \hhh{4mm}
\def \vvv{5mm}
\draw[-] (\hhh*0,\vvv*0) rectangle (\hhh*1,\vvv*1);
\draw[-] (\hhh*1,\vvv*0) rectangle (\hhh*2,\vvv*1);
\draw[-] (\hhh*2,\vvv*0) rectangle (\hhh*3,\vvv*1);
\draw[-] (\hhh*3,\vvv*0) rectangle (\hhh*4,\vvv*1);
\draw[-] (\hhh*1,-\vvv*1) rectangle (\hhh*2,\vvv*0);
\draw[-] (\hhh*2,-\vvv*1) rectangle (\hhh*3,\vvv*0);
\draw[-] (\hhh*3,-\vvv*1) rectangle (\hhh*4,\vvv*0);
\draw[-] (\hhh*2,-\vvv*2) rectangle (\hhh*3,-\vvv*1);
\node at (\hhh*0.5,\vvv*0.5) {\small $1$};
\node at (\hhh*1.5,\vvv*0.5) {\small $2$};
\node at (\hhh*2.5,\vvv*0.5) {\small $2$};
\node at (\hhh*3.5,\vvv*0.5) {\small $4$};
\node at (\hhh*1.5,-\vvv*0.5) {\small $3$};
\node at (\hhh*2.5,-\vvv*0.5) {\small $3$};
\node at (\hhh*3.5,-\vvv*0.5) {\small $5$};
\node at (\hhh*2.5,-\vvv*1.5) {\small $4$};
\node[left] at (-\hhh*0,-\vvv*0.2) {$Q_{T_1}=$};
\end{tikzpicture} 
\hskip 3mm
\begin{tikzpicture}[baseline=4mm]
\def \hhh{4mm}
\def \vvv{5mm}
\draw[-] (\hhh*0,\vvv*0) rectangle (\hhh*1,\vvv*1);
\draw[-] (\hhh*1,\vvv*0) rectangle (\hhh*2,\vvv*1);
\draw[-] (\hhh*2,\vvv*0) rectangle (\hhh*3,\vvv*1);
\draw[-] (\hhh*3,\vvv*0) rectangle (\hhh*4,\vvv*1);
\draw[-] (\hhh*1,-\vvv*1) rectangle (\hhh*2,\vvv*0);
\draw[-] (\hhh*2,-\vvv*1) rectangle (\hhh*3,\vvv*0);
\draw[-] (\hhh*3,-\vvv*1) rectangle (\hhh*4,\vvv*0);
\draw[-] (\hhh*2,-\vvv*2) rectangle (\hhh*3,-\vvv*1);
\node at (\hhh*0.5,\vvv*0.5) {\small $1$};
\node at (\hhh*1.5,\vvv*0.5) {\small $2$};
\node at (\hhh*2.5,\vvv*0.5) {\small $2$};
\node at (\hhh*3.5,\vvv*0.5) {\small $3$};
\node at (\hhh*1.5,-\vvv*0.5) {\small $3$};
\node at (\hhh*2.5,-\vvv*0.5) {\small $4$};
\node at (\hhh*3.5,-\vvv*0.5) {\small $4$};
\node at (\hhh*2.5,-\vvv*1.5) {\small $5$};
\node[left] at (-\hhh*0,-\vvv*0.2) {$Q_{T_2}=$};
\end{tikzpicture} 
\end{displaymath}
\vskip 2mm
\noindent
Moreover, we have
\begin{displaymath}
s_{(5,5,4,3,1)/\delta_4} = 2 P_{(4,3,1)} + P_{(5,2,1)} + P_{(5,3)}.
\end{displaymath}
\end{exm}
\vskip 3mm

\subsection{Ardila-Serrano's expansion of $s_{\delta_{r+1}/\mu}$}
We fix a non-negative integer $r$.
For $\mu\in\cP$ with $D_\mu \subseteq D_{\delta_{r+1}}$,
let us recall the result on the Schur $P$-expansion of he skew Schur function $s_{\delta_{r+1}/\mu}$ by Ardila and Serrano \cite{AS}. 

Let $N=|\delta_{r+1}|-|\mu|$, and 
let $T_{\delta_{r+1}/\mu}$ be the tableau obtained by filling $\delta_{r+1}/\mu$ with $1,2,\ldots,N$ subsequently, starting from the bottom row to top, and from left to right in each row.
For instance,
\begin{displaymath}
\begin{tikzpicture}[baseline=0mm]
\def \hhh{3.8mm}    
\def \vvv{4.5mm}    
\node[left] at (-\hhh*0.2,-\vvv*1) {$T_{\delta_5/(4,1,1)}=$};
\draw[-,black!20,fill=black!10] (\hhh*0,0) rectangle (\hhh*1,\vvv*1);
\draw[-,black!20,fill=black!10] (\hhh*1,0) rectangle (\hhh*2,\vvv*1);
\draw[-,black!20,fill=black!10] (\hhh*2,0) rectangle (\hhh*3,\vvv*1);
\draw[-,black!20,fill=black!10] (\hhh*3,0) rectangle (\hhh*4,\vvv*1);
\draw[-] (\hhh*4,0) rectangle (\hhh*5,\vvv*1);
\draw[-,black!20,fill=black!10] (\hhh*0,-\vvv*1) rectangle (\hhh*1,\vvv*0);
\draw[-] (\hhh*1,-\vvv*1) rectangle (\hhh*2,\vvv*0);
\draw[-] (\hhh*2,-\vvv*1) rectangle (\hhh*3,\vvv*0);
\draw[-] (\hhh*3,-\vvv*1) rectangle (\hhh*4,\vvv*0);
\draw[-,black!20,fill=black!10] (\hhh*0,-\vvv*2) rectangle (\hhh*1,-\vvv*1);
\draw[-] (\hhh*1,-\vvv*2) rectangle (\hhh*2,-\vvv*1);
\draw[-] (\hhh*2,-\vvv*2) rectangle (\hhh*3,-\vvv*1);
\draw[-] (\hhh*0,-\vvv*3) rectangle (\hhh*1,-\vvv*2);
\draw[-] (\hhh*1,-\vvv*3) rectangle (\hhh*2,-\vvv*2);
\draw[-] (\hhh*0,-\vvv*4) rectangle (\hhh*1,-\vvv*3);
\node at (\hhh*0.5,-\vvv*3.5) {\small $1$};
\node at (\hhh*0.5,-\vvv*2.5) {\small $2$};
\node at (\hhh*1.5,-\vvv*2.5) {\small $3$};
\node at (\hhh*1.5,-\vvv*1.5) {\small $4$};
\node at (\hhh*2.5,-\vvv*1.5) {\small $5$};
\node at (\hhh*1.5,-\vvv*0.5) {\small $6$};
\node at (\hhh*2.5,-\vvv*0.5) {\small $7$};
\node at (\hhh*3.5,-\vvv*0.5) {\small $8$};
\node at (\hhh*4.5,\vvv*0.5) {\small $9$};
\node at (\hhh*5.5,-\vvv*4) {.};
\end{tikzpicture}
\end{displaymath}

For $\nu\in\cP^+$ with $|\nu|=N$, let ${\tt B}_{\delta_{r+1}/\mu \, \nu}$ be the set of tableaux $Q$  such that
\begin{itemize}
\item[(1)] $Q\in SST^+_{[N]}(\nu)$ where each entry $i\in [N]$ occurs exactly once,

\item[(2)] if $j$ is directly above $i$ in $T_{\delta_{r+1}/\mu}$, then
$j$ is placed strictly to the right of $i$ in $Q$, 

\item[(3)] if $i+1$ is placed to the right of $i$ in $T_{\delta_{r+1}/\mu}$, then
$i+1$ is strictly below  $i$ in $Q$. 
\end{itemize}
 
\begin{thm}\label{thm:AS}{\rm (\cite[Theorem 4.10]{AS})}
For $\mu\in\cP$ with $D_\mu \subseteq D_{\delta_{r+1}}$,
the skew Schur function $s_{\delta_{r+1}/\mu}$ is given by a non-negative integral 
linear combination of Schur $P$-functions
$$
s_{\delta_{r+1}/\mu} = \sum_{\nu \in \cP^+} b_{\delta_{r+1}/\mu \, \nu} P_\nu,
$$
where $b_{\delta_{r+1}/\mu \, \nu} = |{\tt B}_{\delta_{r+1}/\mu \, \nu}|$.
\end{thm}

Now we show that Theorem \ref{Thm:Main-sect5} (after a little modification of its proof) implies Theorem \ref{thm:AS}.
Let $\la \in \cP$ be such that $D_{\delta_r} \subseteq D_\la \subseteq D_{((r+1)^{r+1})}$.

Let $\nu\in\cP^+$ with $|\nu|=N=|\la|-|\delta_r|$, and let ${\tt L}_{\la/\delta_r\, \nu}$ be as in \eqref{eq:L_{la del}}. Then $|{\tt L}_{\la/\delta_r\, \nu}|=a_{\la/\delta_r\,\nu}$ by \eqref{eq:decomp of stair}.
Let $T \in {\tt L}_{\la/\delta_r\, \nu}$ be given with $w(T)=u_1\cdots u_N$. 
Recall by Lemma \ref{lem:q(n)-lattice property} that  there exists a sequence of strict partitions $\nu^{(m)}\in \cP_n^+$  for $1\leq m\leq N$ such that $u_{N-m+1} \cdots u_N \equiv L^{\nu^{(m)}}$, where $\nu^{(1)}=(1)$, $\nu^{(N)} = \nu$, and
$\nu^{(m)}$ is obtained by adding a box in the $(n-u_m+1)$-st row of $\nu^{(m-1)}$ with $\nu^{(0)}=\emptyset$.

We define $Q'_T$ to be the tableau of shifted shape $\nu$ such that $\nu^{(m)}/\nu^{(m-1)}$ is filled with $m$ for $1\leq m\leq N$. Then we have the following.
 
\begin{thm}\label{Thm:second-sect5}
Let $\la \in \cP$ be such that $D_{\delta_r} \subseteq D_\la \subseteq D_{((r+1)^{r+1})}$ 
and $\nu\in\cP^+$. 
Then we have a bijection 
\begin{equation*}
\xymatrixcolsep{2pc}\xymatrixrowsep{0pc}\xymatrix{
\makebox[3em][c]{${\tt L}_{\la/\delta_r\, \nu}$} \ \ar[r] \ar@{->}[r]  &  \ 
\makebox[4em][l]{${\tt B}_{\delta_{r+1}/(\la^\comp)'\, \nu}$}  \\
\makebox[3em][c]{$T$}    \ \ar[r] \ar@{|->}[r] &  \  \makebox[4em][c]{$Q'_T$}}
\end{equation*}
where 
$\la^\comp:= (r+1-\la_{r+1},r+1-\la_{r},\ldots,r+1-\la_1)$ 
is the complement of $\la$ in $((r+1)^{r+1})$.
\end{thm}
\pf Let $T'_{\la/\delta_r}$ be the tableau obtained by filling $\la/\delta_r$ with $1,2,\ldots,N$ subsequently, starting from the leftmost column to rightmost, and from bottom to top in each column. For instance, when $\la=(5,4,4,4,2)$ and $r=4$, we have 
\begin{displaymath}
\begin{tikzpicture}[baseline=0mm]
\def \hhh{3.8mm}    
\def \vvv{4.5mm}    
\node[left] at (-\hhh*0.2,-\vvv*1) {$T'_{\la/\delta_r}=$};
\draw[-,black!20,fill=black!10] (\hhh*0,0) rectangle (\hhh*1,\vvv*1);
\draw[-,black!20,fill=black!10] (\hhh*1,0) rectangle (\hhh*2,\vvv*1);
\draw[-,black!20,fill=black!10] (\hhh*2,0) rectangle (\hhh*3,\vvv*1);
\draw[-,black!20,fill=black!10] (\hhh*3,0) rectangle (\hhh*4,\vvv*1);
\draw[-] (\hhh*4,0) rectangle (\hhh*5,\vvv*1);
\draw[-,black!20,fill=black!10] (\hhh*0,-\vvv*1) rectangle (\hhh*1,\vvv*0);
\draw[-,black!20,fill=black!10] (\hhh*1,-\vvv*1) rectangle (\hhh*2,\vvv*0);
\draw[-,black!20,fill=black!10] (\hhh*2,-\vvv*1) rectangle (\hhh*3,\vvv*0);
\draw[-] (\hhh*3,-\vvv*1) rectangle (\hhh*4,\vvv*0);
\draw[-,black!20,fill=black!10] (\hhh*0,-\vvv*2) rectangle (\hhh*1,-\vvv*1);
\draw[-,black!20,fill=black!10] (\hhh*1,-\vvv*2) rectangle (\hhh*2,-\vvv*1);
\draw[-] (\hhh*2,-\vvv*2) rectangle (\hhh*3,-\vvv*1);
\draw[-] (\hhh*3,-\vvv*2) rectangle (\hhh*4,-\vvv*1);
\draw[-,black!20,fill=black!10] (\hhh*0,-\vvv*3) rectangle (\hhh*1,-\vvv*2);
\draw[-] (\hhh*1,-\vvv*3) rectangle (\hhh*2,-\vvv*2);
\draw[-] (\hhh*2,-\vvv*3) rectangle (\hhh*3,-\vvv*2);
\draw[-] (\hhh*3,-\vvv*3) rectangle (\hhh*4,-\vvv*2);
\draw[-] (\hhh*0,-\vvv*4) rectangle (\hhh*1,-\vvv*3);
\draw[-] (\hhh*1,-\vvv*4) rectangle (\hhh*2,-\vvv*3);
\node at (\hhh*4.5,\vvv*0.5) {\small $9$};
\node at (\hhh*3.5,-\vvv*0.5) {\small $8$};
\node at (\hhh*3.5,-\vvv*1.5) {\small $7$};
\node at (\hhh*2.5,-\vvv*1.5) {\small $5$};
\node at (\hhh*3.5,-\vvv*2.5) {\small $6$};
\node at (\hhh*2.5,-\vvv*2.5) {\small $4$};
\node at (\hhh*1.5,-\vvv*2.5) {\small $3$};
\node at (\hhh*1.5,-\vvv*3.5) {\small $2$};
\node at (\hhh*0.5,-\vvv*3.5) {\small $1$};
\node at (\hhh*5.5,-\vvv*4) {.};
\end{tikzpicture}
\end{displaymath}
By definition of $Q'_T$, we can check that
\begin{itemize}
\item[(1)] $Q'_T\in SST^+_{[N]}(\nu)$ where each entry $i\in [N]$ occurs exactly once,

\item[(2)] if $j$ is directly above $i$ in $T'_{\la/\delta_r}$, then
then $j$ is strictly below  $i$ in $Q'_T$, 

\item[(3)] if $i+1$ is placed to the right of $i$ in $T'_{\la/\delta_r}$, then
$i+1$ is placed strictly to the right of $i$ in $Q'_T$. 
\end{itemize}

We see that $T_{\delta_{r+1}/(\la^\comp)'}$ is obtained from $T'_{\la/\delta_r}$ by flipping with respect to the main anti-diagonal. This implies that $Q'_T\in {\tt B}_{\delta_{r+1}/(\la^\comp)'\, \nu}$. 
Since the correspondence $T \mapsto Q'_T$ is reversible, it is a bijection.
\qed
\vskip 2mm


\begin{cor}\label{cor:A=B}
Under the above hypothesis, we have a bijection
\begin{equation*}
\xymatrixcolsep{2pc}\xymatrixrowsep{0pc}\xymatrix{
\makebox[3em][c]{${\tt A}_{\la/\delta_r \, \nu}$} \ \ar[r] \ar@{->}[r]  & \ 
\makebox[5em][l]{${\tt B}_{\delta_{r+1}/(\la^\comp)'\, \nu}$}  \\
\makebox[3em][c]{$Q_T$}    \ \ar[r] \ar@{|->}[r] &    \makebox[5em][c]{$Q'_T$}}
\end{equation*}
for $T\in {\tt L}_{\la/\delta_r\, \nu}$.
\end{cor}

Recall that for a skew shape $\eta/\zeta$, we have $s_{\eta/\zeta} = s_{(\eta/\zeta)^\rot}$,
where $(\eta/\zeta)^\rot$ is the (skew) diagram obtained from $\eta/\zeta$
by rotating 180 degree (which can be seen for example by reversing the linear ordering on $\mathbb{N}$ in~\cite{BSS}). Also if $s_{\eta/\zeta}$ has a Schur $P$-expansion, then we have $s_{\eta/\zeta}=s_{\eta'/\zeta'}$ by applying the involution $\omega$ on the ring symmetric function sending $s_\eta$ to $s_{\eta'}$ since $\omega(P_\nu) = P_{\nu}$  for $\nu\in \cP^+$(see~\cite[p. 259, Exercise 3.(a)]{Mac}). 

Hence we have 
\begin{equation*}
s_{\la/\delta_r} = s_{\delta_{r+1}/\la^\comp}=s_{\delta_{r+1}/(\la^\comp)'},
\end{equation*} 
for $\la \in \cP$ such that $D_{\delta_r} \subseteq D_\la \subseteq D_{((r+1)^{r+1})}$.
This implies that
\begin{equation}\label{eq:a=b}
a_{\la / \delta_{r} \, \nu} 
= b_{\delta_{r+1}/\la^\comp\, \nu}
= b_{\delta_{r+1}/(\la^\comp)'\, \nu},
\end{equation}
for $\nu\in\cP^+$, where $a_{\la / \delta_{r} \, \nu}$ are given in \eqref{Eq:s-to-P}. Equivalently, we have   
\begin{equation}\label{eq:a=b-2}
a_{(\mu^\comp)' / \delta_{r} \, \nu} 
= b_{\delta_{r+1}/\mu'\, \nu}
= b_{\delta_{r+1}/\mu\, \nu},
\end{equation}
for $\mu\in\cP$ with $D_\mu \subseteq D_{\delta_{r+1}}$. Therefore Theorem \ref{thm:AS} follows from Theorem \ref{Thm:Main-sect5}, Corollary \ref{cor:A=B}, and \eqref{eq:a=b} (or \eqref{eq:a=b-2}).

\section{Schur expansion of Schur $P$-function}

For $\la\in\cP^+$ and $\mu\in \cP$, let $g_{\lambda \mu}$ be the coefficient of $s_\mu$ in the Schur expansion of $P_\la$, that is,
\begin{equation}
P_\lambda = \sum_{\mu} g_{\lambda \mu} s_\mu.
\end{equation}
The purpose of this section is to give an alternate proof of the following combinatorial description of $g_{\la\mu}$ due to Stembridge.  

\begin{thm}\label{thm-coeff-g}
{\rm (\cite[Theorem 9.3]{Ste})}
For $\la \in \cP^+$ and $\mu \in \cP$, we have
\begin{equation*}
g_{\la\mu} = \left| {\tt G}_{\la\mu} \right|,
\end{equation*}
where ${\tt G}_{\la\mu}$ is the set of tableaux $Q$ such that
\begin{itemize}
\item[\rm (1)] $Q\in SST_{\mc N}(\mu)$ with $c_{k}(Q)+c_{k'}(Q)=\la_k$ for $k\geq 1$, 

\item[\rm (2)] for $k\geq 1$, if $x$ is the rightmost letter in $w(Q)$ with $|x|=k$, then $x=k$,

\item[\rm (3)] $w(Q)$ satisfies the lattice property.
\end{itemize}
\end{thm}

\pf The proof is similar to that of Theorem \ref{Thm:our-characterization}.
Choose $n$ such that $\la\in \cP_n^+$ and $\mu \in \cP_n$.
Let
\begin{equation*}
{\tt \bf L}_{\la\mu}=
\{\,T\,|\,T\in \B_n(\la),\,\tf_iT={\bf 0}\ (1\leq i\leq n-1),\ {\rm wt} (T) =w_0\mu \,\}.
\end{equation*}
Then we have as a $\gl(n)$-crystal
\begin{equation}
\B_n(\la)\cong \bigsqcup_{\mu}B_n(\mu)^{\oplus |{\tt \bf L}_{\la\mu}|},
\end{equation}
and hence $g_{\la\mu}=|{\tt \bf L}_{\lambda \mu}|$ by linear independence of Schur polynomials.
Let us define a map
\begin{equation*}\label{eq:LRS2 recording}
\xymatrixcolsep{2pc}\xymatrixrowsep{2mm}\xymatrix@1{
\makebox[3em][c]{${\tt \bf L}_{\lambda \mu}$}      \ar[r] \ar@{->}[r]  &   \makebox[3em][c]{${\tt G}_{\la \mu}$}  \\
\makebox[3em][c]{$T$}    \ar[r] \ar@{|->}[r] & \    \makebox[3em][c]{$Q_T$}}
\end{equation*}
as follows. Let $T \in {\tt \bf L}_{\la \mu}$ be given.
Assume that $w_\rev(T)=u_1\cdots u_N$ where $N=|\lambda|$. 
Since $T$ is a $\gl(n)$-lowest weight vector, we have by \eqref{eq:tensor product of crystals} that $u_{N-m+1}\otimes  \cdots \otimes u_N \in \B^{\otimes m}_n$ is a $\gl(n)$-lowest weight element for $1\leq m\leq N$. 
This implies that there exists $\mu^{(m)}\in \cP_n$ for $1\leq m\leq N$ such that 
$u_{N-m+1} \cdots u_N$ is equivalent as an element of $\gl(n)$-crystal to a $\gl(n)$-lowest weight element in $B_n(\mu^{(m)})$,
where $\mu^{(N)} = \mu$ and
$\mu^{(m)}$ is obtained by adding a box in the $(n-u_m+1)$-st row of $\mu^{(m-1)}$ with $\mu^{(0)}=\emptyset$. 

We define $Q_T$ to be a tableau of shape $\mu$ with entries in $\mc N$, where $\mu^{(m)}/\mu^{(m-1)}$ is filled with
\begin{equation*}\label{eq:labeling of Q_T-2}
\begin{cases}
k', & \text{if $u_m$ belongs to $T^{(k)}\!\!\uparrow$},\\
k, & \text{if $u_m$ belongs to $T^{(k)}\!\!\downarrow$},\\
\end{cases}
\end{equation*}
for some $1\leq k\leq \ell(\la)$. By almost the same arguments as in the proof of Theorem \ref{Thm:our-characterization}, we see that $Q_T$ satisfies the conditions (1) and (2) for $\tt{G}_{\la\mu}$, and $w(Q_T)$ satisfies the ``lattice property", which implies that it satisfies the lattice property by Corollary \ref{cor:equivalence of lattice conditions}. (We leave the details to the reader.) Finally the correspondence $T \mapsto Q_T$ is a well-defined bijection. \qed

\begin{exm}
Let $\la = (3,1)$.
From Figure~\ref{Fig:q(3)-crystal}
we get three $\gl(3)$-lowest weight vectors in $\B_3(\la)$
\begin{center}
\begin{tikzpicture}[scale=0.8]
\def \hhh{35mm}
\def \vvv{0}
\def \widthh{.5ex}
\def \highhh{-.3ex}
\node (node_56) at (\hhh*2,-\vvv*0) [draw,draw=none] {${\def\lr#1{\multicolumn{1}{|@{\hspace{\widthh}}c@{\hspace{\widthh}}|}{\raisebox{\highhh}{$#1$}}}\raisebox{\highhh}
{$\begin{array}[b]{*{3}c}\cline{1-3}\lr{3}&\lr{2}&\lr{3}\\\cline{1-3}&\lr{1}\\\cline{2-2}\end{array}$}}$.};
\node (node_64) at (\hhh,-\vvv*0) [draw,draw=none] {${\def\lr#1{\multicolumn{1}{|@{\hspace{\widthh}}c@{\hspace{\widthh}}|}{\raisebox{\highhh}{$#1$}}}\raisebox{\highhh}
{$\begin{array}[b]{*{3}c}\cline{1-3}\lr{3}&\lr{2}&\lr{3}\\\cline{1-3}&\lr{2}\\\cline{2-2}\end{array}$}}$};
\node (node_73) at (0,-\vvv*0) [draw,draw=none] {${\def\lr#1{\multicolumn{1}{|@{\hspace{\widthh}}c@{\hspace{\widthh}}|}{\raisebox{\highhh}{$#1$}}}\raisebox{\highhh}
{$\begin{array}[b]{*{3}c}\cline{1-3}\lr{3}&\lr{3}&\lr{3}\\\cline{1-3}&\lr{2}\\\cline{2-2}\end{array}$}}$};
\end{tikzpicture}
\end{center}
By applying the mapping $T \mapsto Q_T$ in the proof of Theorem~\ref{thm-coeff-g} to these tableaux we have 
\begin{center}
\begin{tikzpicture}[scale=0.8,baseline=0mm]
\def \hhh{35mm}
\def \vvv{2ex}
\def \widthh{.5ex}
\def \highhh{-.3ex}
\node (node_56) at (\hhh*2,-\vvv*1) [draw,draw=none] {${\def\lr#1{\multicolumn{1}{|@{\hspace{\widthh}}c@{\hspace{\widthh}}|}{\raisebox{\highhh}{$#1$}}}\raisebox{\highhh}
{$\begin{array}[b]{*{2}c}
\cline{1-2}\lr{1'}&\lr{\hskip 0.5mm 1 \hskip 0.5mm}\\
\cline{1-2}\lr{1}&\\
\cline{1-1}\lr{2}&\\
\cline{1-1}
\end{array}$}}$};
\node (node_64) at (\hhh,-\vvv*0) [draw,draw=none] {${\def\lr#1{\multicolumn{1}{|@{\hspace{\widthh}}c@{\hspace{\widthh}}|}{\raisebox{\highhh}{$#1$}}}\raisebox{\highhh}
{$\begin{array}[b]{*{2}c}
\cline{1-2}\lr{1'}&\lr{\hskip 0.5mm 1\hskip 0.5mm }\\
\cline{1-2}\lr{1}&\lr{2}\\
\cline{1-2}
\end{array}$}}$};
\node (node_73) at (0,-\vvv*0) [draw,draw=none] {${\def\lr#1{\multicolumn{1}{|@{\hspace{\widthh}}c@{\hspace{\widthh}}|}{\raisebox{\highhh}{$#1$}}}\raisebox{\highhh}
{$\begin{array}[b]{*{3}c}\cline{1-3}\lr{1}&\lr{1}&\lr{1}\\\cline{1-3}\lr{2}& &\\\cline{1-1}\end{array}$}}$};
\end{tikzpicture}.
\end{center}
Thus $P_{(3,1)} = s_{(3,1)} + s_{(2,2)} + s_{(2,1,1)}.$
\end{exm}

\begin{rem}
Let $\la \in \cP^+$ be such that $D^+_\la \subseteq D^+_{\delta_{r+1}}$ for some $r \geq 0$. Let $\la^{\compp}$ be a strict partition obtained by counting complementary boxes $D^+_{\delta_{r+1}}\setminus D^+_\la$ in each column from right to left.
It is shown in \cite{Dew} that 
\begin{equation*}
s_{\delta_{r+1}/\la} = \sum_{\substack{\ \nu \in \cP^+ \\ |\nu|=|\la|}} g_{\nu \la} P_{\nu^{\compp}}.
\end{equation*} 
By \eqref{eq:a=b} or \eqref{eq:a=b-2}, we have 
$g_{\nu  \,  \la}=a_{\la^\comp/\delta_r\, (\nu^\compp)'}.$
One may expect that there is a natural bijection between
${\tt G}_{\nu\, \la}$ and ${\tt A}_{\la^\comp/\delta_r\, (\nu^\compp)'}$, but we do not know the answer yet.
\end{rem}

\section{Semistandard decomposition tableaux of skew shapes}
\label{Sec:skew-SSDT}
Let $\la/\mu$ be a shifted skew diagram for $\la, \mu \in \cP^+$ with $D^+_\mu \subseteq D^+_\la$. Without loss of generality, we assume in this section that $\la_1 > \mu_1$ and $\ell(\la) > \ell(\mu)$.

Let $T$ be a tableau of shifted skew shape $\la/\mu$.
For $p,q\geq 1$, let $T(p,q)$ denote the entry of $T$ at the $p$-th row and the $q$-th diagonal from the main diagonal in $D^+_\la$ (that is, $\{\,(i,i)\,|\,i\geq 1 \}\cap  D^+_\la$) whenever it is defined. Note that $T(p,q)$ is not necessarily equal to $T_{p,q}$ if $\mu$ is nonempty. 

For example, when $\la/\mu= (5,4,2) / (3,1)$, we have
\vskip 5mm

\begin{center}
\begin{tikzpicture}[baseline=0mm]
\def \hhh{9.5mm}
\def \vvv{6mm}
\draw[-,black!20,fill=black!10] (\hhh*0,\vvv*0) rectangle (\hhh*1,\vvv*1) node[pos=0.5] {\scalebox{.7}{ }};
\draw[-,black!20,fill=black!10] (\hhh*1,\vvv*0) rectangle (\hhh*2,\vvv*1) node[pos=0.5] {\scalebox{.7}{ }};
\draw[-,black!20,fill=black!10] (\hhh*2,\vvv*0) rectangle (\hhh*3,\vvv*1) node[pos=0.5] {\scalebox{.7}{ }};
\draw (\hhh*3,\vvv*0) rectangle (\hhh*4,\vvv*1) node[pos=0.5] {\scalebox{.7}{$\tiny T(1,4)$}};
\draw (\hhh*4,\vvv*0) rectangle (\hhh*5,\vvv*1) node[pos=0.5] {\scalebox{.7}{$T(1,5)$}};
\draw[-,black!20,fill=black!10] (\hhh*1,-\vvv*1) rectangle (\hhh*2,\vvv*0) node[pos=0.5] {\scalebox{.7}{ }};
\draw (\hhh*2,-\vvv*1) rectangle (\hhh*3,\vvv*0) node[pos=0.5] {\scalebox{.7}{$T(2,2)$}};
\draw (\hhh*3,-\vvv*1) rectangle (\hhh*4,\vvv*0) node[pos=0.5] {\scalebox{.7}{$T(2,3)$}};
\draw (\hhh*4,-\vvv*1) rectangle (\hhh*5,\vvv*0) node[pos=0.5] {\scalebox{.7}{$T(2,4)$}};
\draw (\hhh*2,-\vvv*2) rectangle (\hhh*3,-\vvv*1) node[pos=0.5] {\scalebox{.7}{$T(3,1)$}};
\draw (\hhh*3,-\vvv*2) rectangle (\hhh*4,-\vvv*1) node[pos=0.5] {\scalebox{.7}{$T(3,2)$}};
\end{tikzpicture}
\hskip 6mm
\begin{tikzpicture}[baseline=0mm]
\def \vvv{6mm}
\node[left] at (0,-\vvv*0.5) {$=$};
\end{tikzpicture}
\hskip 6mm
\begin{tikzpicture}[baseline=0mm]
\def \hhh{9mm}
\def \vvv{6mm}
\draw[-,black!20,fill=black!10] (\hhh*0,\vvv*0) rectangle (\hhh*1,\vvv*1) node[pos=0.5] {\scalebox{.9}{ }};
\draw[-,black!20,fill=black!10] (\hhh*1,\vvv*0) rectangle (\hhh*2,\vvv*1) node[pos=0.5] {\scalebox{.9}{ }};
\draw[-,black!20,fill=black!10] (\hhh*2,\vvv*0) rectangle (\hhh*3,\vvv*1) node[pos=0.5] {\scalebox{.9}{ }};
\draw (\hhh*3,\vvv*0) rectangle (\hhh*4,\vvv*1) node[pos=0.5] {\scalebox{.8}{$T_{1,1}$}};
\draw (\hhh*4,\vvv*0) rectangle (\hhh*5,\vvv*1) node[pos=0.5] {\scalebox{.8}{$T_{1,2}$}};
\draw[-,black!20,fill=black!10] (\hhh*1,-\vvv*1) rectangle (\hhh*2,\vvv*0) node[pos=0.5] {\scalebox{.9}{ }};
\draw (\hhh*2,-\vvv*1) rectangle (\hhh*3,\vvv*0) node[pos=0.5] {\scalebox{.8}{$T_{2,1}$}};
\draw (\hhh*3,-\vvv*1) rectangle (\hhh*4,\vvv*0) node[pos=0.5] {\scalebox{.8}{$T_{2,2}$}};
\draw (\hhh*4,-\vvv*1) rectangle (\hhh*5,\vvv*0) node[pos=0.5] {\scalebox{.8}{$T_{2,3}$}};
\draw (\hhh*2,-\vvv*2) rectangle (\hhh*3,-\vvv*1) node[pos=0.5] {\scalebox{.8}{$T_{3,1}$}};
\draw (\hhh*3,-\vvv*2) rectangle (\hhh*4,-\vvv*1) node[pos=0.5] {\scalebox{.8}{$T_{3,2}$}};
\end{tikzpicture}
\end{center}
\vskip 3mm

\begin{df}\label{def:skew SSDT}
For $\la, \mu \in \cP^+$ with $D^+_\mu \subseteq D^+_\la$,
a {\em skew semistandard decomposition tableau} $T$ of shape $\la/\mu$
is a tableau of shifted shape $\la/\mu$ with entries in $\mathbb{N}$ such that 
$T^{(k)}$ is a hook word for $1\leq k\leq \ell(\la)$ and 
the following holds for $1\leq k < \ell(\la)$ and $1 \leq i \leq j\leq  \la_{k+1}$:
\begin{itemize}
\item[(S1)] if $T({k,i}) \leq T({k+1,j})$, 
then $i\neq 1$ and $T({k+1,i-1}) < T({k+1,j})$, 

\item[(S2)] if $T({k,i}) > T({k+1,j})$, then $T({k,i}) \geq T({k,j+1})$,
\end{itemize} 
where we assume that $T(p,q)$ for $p,q\geq 1$ is empty if it is not defined. 
\end{df} 

Let $SSDT(\la/\mu)$ be the set consisting of skew semistandard decomposition tableaux of shape $\la/\mu$. Note that when $\mu$ is empty, the set $SSDT(\la/\mu)$ is equal to $SSDT(\la)$ by Proposition~\ref{prop:SSDT-checker}.

Suppose that $\ell(\la)\leq n$.
Let $\B_n(\la/\mu)$ be the set of $T\in SSDT(\la/\mu)$ with entries in $[n]$. As in \eqref{eq:embedding of SSDT}, consider the injective map
\begin{equation}\label{eq:embedding of skew SSDT}
\xymatrixcolsep{2pc}\xymatrixrowsep{0pc}\xymatrix{
\makebox[4em][c]{$\B_n(\la/\mu)$} \  \ar@{^{(}->}[r]  &\  \makebox[3em][c]{$\W_{[n]}$ \ }  \\
 \makebox[4em][c]{$T$}  \ar@{|->}[r] & \   \makebox[3em][l]{$w_\rev(T)$.}
}
\end{equation}

\begin{prop}\label{prop:skew SSDT}
Under the above hypothesis, the image of $\B_n(\lambda/\mu)$ in \eqref{eq:embedding of skew SSDT} together with $\{{\bf 0}\}$ is invariant under the action of $\te_i$ and $\tf_i$ for $i\in I$, and hence $\B_n(\la/\mu)$ is a $\q(n)$-crystal.
\end{prop}
\pf Choose a sufficiently large $M$ such that all the entries in $L^\mu_M$ are greater than $n$. For a tableau $T$ of shifted shape $\la/\mu$ with entries in $[n]$, let $\widetilde{T}:=L^\mu_M\ast T$ be the tableau of shifted shape $\la$, that is, the subtableau of shape shifted $\mu$ in $\widetilde{T}$ is $L^\mu_M$ and its complement in $\widetilde{T}$ is $T$. By definition of $SSDT(\la/\mu)$ and Proposition \ref{prop:SSDT-checker}, we have
\begin{equation}\label{eq:skew checker}
\text{$T\in \B_n(\la/\mu)$ \ if and only if \ $\widetilde{T}\in \B_M(\la)$.}
\end{equation}

Let $T\in \B_n(\la/\mu)$ and $i\in I$ be given. If $\tilde{x}_i \widetilde{T}\neq {\bf 0}$ ($x=e, f$), then we have by \eqref{eq:skew checker} that $\tilde{x}_i \widetilde{T} = L^\mu_M\ast T'$ for some $T'\in \B_n(\la/\mu)$. This implies that $\tilde{x}_i w_{\rm rev}(T)=w_{\rm rev}(T')$. Therefore, the image of $\B_n(\lambda/\mu)$ in \eqref{eq:embedding of skew SSDT} together with $\{{\bf 0}\}$ is invariant under the action of $\te_i$ and $\tf_i$ for $i\in I$.
\qed
\vskip 3mm

Since $\B_n(\la/\mu)$ is a subcrystal of $\B_n^{\otimes N}$ with $N=|\la|-|\mu|$, we have  
\begin{equation}\label{eq:skew decomp}
\B_n(\la/\mu) \cong \bigsqcup_{\substack{\nu\in\cP_n^+ \\ |\nu|=N}}\B_n(\nu)^{\oplus f^{\la/\mu}_{\nu}(n)}
\end{equation}
for some $f^{\la/\mu}_{\nu}(n)\in \Z_+$. Moreover by Remark \ref{rem:stability}, we have
\begin{equation}\label{eq:skew decomp mult}
f^{\la/\mu}_{\nu}:=f^{\la/\mu}_{\nu}(m)=f^{\la/\mu}_{\nu}(n)\quad (m\geq n).
\end{equation}
If we put
\begin{equation*}
{P}^\circ_{\la/\mu}=\sum_{T\in SSDT(\la/\mu)}x^T,
\end{equation*}
then we have from \eqref{eq:skew decomp} and \eqref{eq:skew decomp mult}
\begin{equation}\label{eq:Schur P-expansion of skew SSDT}
{P}^\circ_{\la/\mu}=\sum_{\nu\in\cP^+}f^{\la/\mu}_{\nu}P_\nu.
\end{equation}

\begin{exm}{\rm
For $\eta\in \cP_n^+$ with $\ell(\eta)=\ell$, let $\la=\eta+L\delta_\ell$ and $\mu=L\delta_{\ell}\in \cP^+_n$, where $L\geq \eta_1$. 
 Since each column in $\la/\mu$ has at most one box, we have
\begin{equation*}
\B_n(\la/\mu)\cong \B_n(\eta_1)\otimes \cdots \otimes \B_n(\eta_\ell). 
\end{equation*}
By applying Theorem~\ref{Thm:our-characterization} repeatedly, we see that $f^{\la/\mu}_{\nu}$ for $\nu\in\cP_n^+$ in this case is equal to the number of tableaux $Q$ such that
\begin{itemize}
\item[(1)] $Q\in SST^+_{\mc N}(\nu)$ with $c_{k}(Q)+c_{k'}(Q)=\eta_k$ for $k\geq 1$,

\item[(2)] for each $k\geq 1$, if $x$ is the rightmost in $w(Q)$ with $|x|=k$, then $x=k$.
\end{itemize}
}
\end{exm}
 
One can generalize the notion of ``lattice property" in Definition \ref{Def:new lattice property} to describe the coefficient $f^{\la/\mu}_{\nu}$.

\begin{df}
\label{def:skew lattice property}
Let $w=w_1\cdots w_N\in \W_{\mc N}$ be given and let $w^*=w_1^*\cdots w_N^*$ be the word associated to $w$ given in Definition \ref{Def:w and w^*}. 
For $\mu\in\cP^+$, we say that $w$ satisfies the ``{\em $\mu$-lattice property}" if $w^*$ satisfies the following for each $k\geq 1$:
\begin{itemize}
\item[(L1)]
if $k>\ell(\mu)$ and $w^*_i=k_1$, then no $k+1_j$ for $j\geq 1$ occurs in 
$w^*_1\cdots w^*_{i-1}$,

\item[(L2)] if $(w^*_s,w^*_t)=(k+1_{i},k_{i+1-\alpha_k})$ for some $s<t$ and $\alpha_k<i$, then no $k+1_j$ for $i<j$ occurs in $w^*_s\cdots w^*_t$, 

\item[(L3)] if $(w^*_s,w^*_t)=(k_{j+1-\alpha_k},k+1_{j})$ for some $s<t$ and $\alpha_k<j$, then no $k_i$ for $i\leq j-\alpha_k$ occurs in $w^*_s\cdots w^*_t$,
\end{itemize}
where $\alpha_k=\mu_k-\mu_{k+1}$.
\end{df}

\begin{thm}\label{thm-our-characterization-skew}
For $\la, \mu, \nu \in \cP^+$, we have 
$$
f^{\la/\mu}_{\nu} = \left|{\tt F}^{\la/\mu}_{\nu}\right|,
$$
where 
${\tt F}^{\la/\mu}_{\nu}$ is the set of tableaux $Q$ such that  
\begin{itemize}
\item[\rm (1)] $Q\in SST^+_{\mc N}(\nu)$ with $c_{k}(Q)+c_{k'}(Q)=\la_k-\mu_k$ for $k\geq 1$,

\item[\rm (2)] for $k\geq 1$, if $x$ is the rightmost letter in $w(Q)$ with $|x|=k$, then $x=k$,

\item[\rm (3)] $w(Q)$ satisfies the $\mu$-lattice property. 

\end{itemize}
\end{thm}
\pf The proof is similar to that of Theorem \ref{Thm:our-characterization}.
Choose $n$ such that $\la, \mu, \nu\in \cP_n^+$. Put
\begin{equation*}
{\tt \bf L}_\nu^{\la/\mu}
=\{\,T\,|\, T \in \B_n(\la/\mu),\ T \equiv L^\nu\,\}.
\end{equation*}
From \eqref{eq:skew decomp} and \eqref{eq:skew decomp mult},  we have $|{\tt \bf L}^{\la/\mu}_{\nu}|=f^{\la/\mu}_{\nu}$. 
Let us define a map
\begin{equation*}\label{eq:LRS1 recording}
\xymatrixcolsep{2pc}\xymatrixrowsep{0pc}\xymatrix{
{\tt \bf L}^{\la/\mu}_{\nu}\ \  \ar@{->}[r]  &\ \ {\tt F}^{\la/\mu}_{\nu}  \\
\ \ \ T \ \ \ \ar@{|->}[r] & \  {Q}_T}
\end{equation*}
as follows. 
Let $N=|\la|-|\mu|$.
Suppose that $T \in {\tt L}^{\la/\mu}_{\nu}$ is given with $w_{\rm rev}(T)=u_1\cdots u_N$. 
By Lemma \ref{lem:q(n)-lattice property} there exists $\nu^{(m)}\in \cP_n^+$ for $1\leq m \leq N$ such that $u_{N-m+1}\cdots u_N \equiv L^{\nu^{(m)}}$
where $\nu^{(N)}=\nu$ and 
$\nu^{(m)}$ is obtained  by adding a box in the $(n-u_m+1)$-st row of $\nu^{(m-1)}$ with $\nu^{(0)}=\emptyset$.

Note that
$w_{\rm rev}(T)=T^{(\ell(\la))}\cdots T^{(1)}$, 
where $T^{(k)}$ is a hook word for $1\leq k\leq \ell(\la)$.
Then we define ${Q}_T$ to be a tableau of shifted shape $\nu$ with entries in $\mc N$, where $\nu^{(m)}/\nu^{(m-1)}$ is filled with
\begin{equation}\label{eq:labeling of ov Q_T}
\begin{cases}
k', & \text{if $u_m$ belongs to $T^{(k)}\!\!\uparrow$},\\
k, & \text{if $u_m$ belongs to $T^{(k)}\!\!\downarrow$},\\
\end{cases}
\end{equation}
for some $1\leq k\leq \ell(\la)$. 

First, by the same argument as in the proof of Theorem \ref{Thm:our-characterization}, we see that $Q_T$ satisfies the condition (2) for $\tt{F}^{\la/\mu}_{\nu}$ by the same argument as in the proof of Theorem \ref{Thm:our-characterization}.   

Let us check that $w({Q}_T)$ satisfies the $\mu$-lattice property. 
If we label $k$ and $k'$ in \eqref{eq:labeling of ov Q_T} 
as $k_j$ and $k'_j$, respectively, when $u_m=T_{k,j}$, then it coincides with the labeling on the letters in $w({Q}_T)$ given in Definition \ref{Def:w and w^*}(1). 

Choose a sufficiently large $M$ such that all the entries in $L^\mu_M$ are greater than $n$. Let $S=L_M^\mu\ast T$ (see the proof of Proposition \ref{prop:skew SSDT}). Since $S\in \B_M(\la)$, the conditions Proposition~\ref{prop:SSDT-checker}(1), (2), and (3) on $S$ and hence on $T$ (cf. \eqref{eq:skew checker}) imply the conditions Definition~\ref{def:skew lattice property}(L1), (L2), and (L3), respectively. 
Therefore, ${Q}_T\in {\tt F}^{\la/\mu}_{\nu}$.

Finally the correspondence $T \mapsto {Q}_T$ is injective and also reversible. Hence it is a bijection. \qed

\begin{rem}
In general, $P^\circ_{\la/\mu}$ is not equal to the usual skew Schur $P$-function $P_{\la/\mu}$, or $f^{\la/\mu}_{\nu}$ is not necessarily equal to $f^\la_{\mu\nu}$. It would be interesting to have a representation-theoretic interpretation of the Schur $P$-expansion of $P^\circ_{\la/\mu}$ \eqref{eq:Schur P-expansion of skew SSDT}.
\end{rem}

\begin{exm}\label{Ex:skew-SSDTs}
Let $\lambda = (6,5,2,1)$ and $\mu=(4,2)$. 
Then 
$$P^\circ_{\lambda/\mu} =2 P_{(6,2)}+3 P_{(5,3)} + 5 P_{(5,2,1)} + 4 P_{(4,3,1)}$$
since $\bigsqcup_{\nu}{\tt F}^{\la/\mu}_{\nu}$ consists of
\vskip 4mm
\noindent
\begin{tikzpicture}[baseline=0mm]
\def \hhh{4mm}
\def \vvv{5mm}
\draw[-] (\hhh*0,\vvv*0) rectangle (\hhh*1,\vvv*1) node[pos=.5] {\small $1$};
\draw[-] (\hhh*1,\vvv*0) rectangle (\hhh*2,\vvv*1) node[pos=.5] {\small $1$};
\draw[-] (\hhh*2,\vvv*0) rectangle (\hhh*3,\vvv*1) node[pos=.5] {\small $2'$};
\draw[-] (\hhh*3,\vvv*0) rectangle (\hhh*4,\vvv*1) node[pos=.5] {\small $2$};
\draw[-] (\hhh*4,\vvv*0) rectangle (\hhh*5,\vvv*1) node[pos=.5] {\small $3$};
\draw[-] (\hhh*5,\vvv*0) rectangle (\hhh*6,\vvv*1) node[pos=.5] {\small $3$};
\draw[-] (\hhh*1,-\vvv*1) rectangle (\hhh*2,\vvv*0) node[pos=.5] {\small $2$};
\draw[-] (\hhh*2,-\vvv*1) rectangle (\hhh*3,\vvv*0) node[pos=.5] {\small $4$};
\end{tikzpicture} 
\hskip 3mm
\begin{tikzpicture}[baseline=0mm]
\def \hhh{4mm}
\def \vvv{5mm}
\draw[-] (\hhh*0,\vvv*0) rectangle (\hhh*1,\vvv*1) node[pos=.5] {\small $1$};
\draw[-] (\hhh*1,\vvv*0) rectangle (\hhh*2,\vvv*1) node[pos=.5] {\small $1$};
\draw[-] (\hhh*2,\vvv*0) rectangle (\hhh*3,\vvv*1) node[pos=.5] {\small $2$};
\draw[-] (\hhh*3,\vvv*0) rectangle (\hhh*4,\vvv*1) node[pos=.5] {\small $2$};
\draw[-] (\hhh*4,\vvv*0) rectangle (\hhh*5,\vvv*1) node[pos=.5] {\small $3$};
\draw[-] (\hhh*5,\vvv*0) rectangle (\hhh*6,\vvv*1) node[pos=.5] {\small $3$};
\draw[-] (\hhh*1,-\vvv*1) rectangle (\hhh*2,\vvv*0) node[pos=.5] {\small $2$};
\draw[-] (\hhh*2,-\vvv*1) rectangle (\hhh*3,\vvv*0) node[pos=.5] {\small $4$};
\end{tikzpicture} 
\hskip 3mm
\begin{tikzpicture}[baseline=0mm]
\def \hhh{4mm}
\def \vvv{5mm}
\node[right] at (\hhh*1,-\vvv*0) {(when $\nu=(6,2)$)};
\end{tikzpicture} 
\vskip 4mm
\noindent
\begin{tikzpicture}[baseline=0mm]
\def \hhh{4mm}
\def \vvv{5mm}
\draw[-] (\hhh*0,\vvv*0) rectangle (\hhh*1,\vvv*1) node[pos=.5] {\small $1$};
\draw[-] (\hhh*1,\vvv*0) rectangle (\hhh*2,\vvv*1) node[pos=.5] {\small $1$};
\draw[-] (\hhh*2,\vvv*0) rectangle (\hhh*3,\vvv*1) node[pos=.5] {\small $2'$};
\draw[-] (\hhh*3,\vvv*0) rectangle (\hhh*4,\vvv*1) node[pos=.5] {\small $3$};
\draw[-] (\hhh*4,\vvv*0) rectangle (\hhh*5,\vvv*1) node[pos=.5] {\small $3$};
\draw[-] (\hhh*1,-\vvv*1) rectangle (\hhh*2,\vvv*0) node[pos=.5] {\small $2$};
\draw[-] (\hhh*2,-\vvv*1) rectangle (\hhh*3,\vvv*0) node[pos=.5] {\small $2$};
\draw[-] (\hhh*3,-\vvv*1) rectangle (\hhh*4,\vvv*0) node[pos=.5] {\small $4$};
\end{tikzpicture} 
\hskip 4mm
\begin{tikzpicture}[baseline=0mm]
\def \hhh{4mm}
\def \vvv{5mm}
\draw[-] (\hhh*0,\vvv*0) rectangle (\hhh*1,\vvv*1) node[pos=.5] {\small $1$};
\draw[-] (\hhh*1,\vvv*0) rectangle (\hhh*2,\vvv*1) node[pos=.5] {\small $1$};
\draw[-] (\hhh*2,\vvv*0) rectangle (\hhh*3,\vvv*1) node[pos=.5] {\small $2'$};
\draw[-] (\hhh*3,\vvv*0) rectangle (\hhh*4,\vvv*1) node[pos=.5] {\small $2$};
\draw[-] (\hhh*4,\vvv*0) rectangle (\hhh*5,\vvv*1) node[pos=.5] {\small $3$};
\draw[-] (\hhh*1,-\vvv*1) rectangle (\hhh*2,\vvv*0) node[pos=.5] {\small $2$};
\draw[-] (\hhh*2,-\vvv*1) rectangle (\hhh*3,\vvv*0) node[pos=.5] {\small $3$};
\draw[-] (\hhh*3,-\vvv*1) rectangle (\hhh*4,\vvv*0) node[pos=.5] {\small $4$};
\end{tikzpicture} 
\hskip 4mm
\begin{tikzpicture}[baseline=0mm]
\def \hhh{4mm}
\def \vvv{5mm}
\draw[-] (\hhh*0,\vvv*0) rectangle (\hhh*1,\vvv*1) node[pos=.5] {\small $1$};
\draw[-] (\hhh*1,\vvv*0) rectangle (\hhh*2,\vvv*1) node[pos=.5] {\small $1$};
\draw[-] (\hhh*2,\vvv*0) rectangle (\hhh*3,\vvv*1) node[pos=.5] {\small $2$};
\draw[-] (\hhh*3,\vvv*0) rectangle (\hhh*4,\vvv*1) node[pos=.5] {\small $2$};
\draw[-] (\hhh*4,\vvv*0) rectangle (\hhh*5,\vvv*1) node[pos=.5] {\small $3$};
\draw[-] (\hhh*1,-\vvv*1) rectangle (\hhh*2,\vvv*0) node[pos=.5] {\small $2$};
\draw[-] (\hhh*2,-\vvv*1) rectangle (\hhh*3,\vvv*0) node[pos=.5] {\small $3$};
\draw[-] (\hhh*3,-\vvv*1) rectangle (\hhh*4,\vvv*0) node[pos=.5] {\small $4$};
\end{tikzpicture} 
\hskip 3mm
\begin{tikzpicture}[baseline=0mm]
\def \hhh{4mm}
\def \vvv{5mm}
\node[right] at (\hhh*1,-\vvv*0) {(when $\nu=(5,3)$)};
\end{tikzpicture} 
\vskip 4mm
\noindent
\begin{tikzpicture}[baseline=0mm]
\def \hhh{4mm}
\def \vvv{5mm}
\draw[-] (\hhh*0,\vvv*0) rectangle (\hhh*1,\vvv*1) node[pos=.5] {\small $1$};
\draw[-] (\hhh*1,\vvv*0) rectangle (\hhh*2,\vvv*1) node[pos=.5] {\small $1$};
\draw[-] (\hhh*2,\vvv*0) rectangle (\hhh*3,\vvv*1) node[pos=.5] {\small $2'$};
\draw[-] (\hhh*3,\vvv*0) rectangle (\hhh*4,\vvv*1) node[pos=.5] {\small $3$};
\draw[-] (\hhh*4,\vvv*0) rectangle (\hhh*5,\vvv*1) node[pos=.5] {\small $3$};
\draw[-] (\hhh*1,-\vvv*1) rectangle (\hhh*2,\vvv*0) node[pos=.5] {\small $2$};
\draw[-] (\hhh*2,-\vvv*1) rectangle (\hhh*3,\vvv*0) node[pos=.5] {\small $2$};
\draw[-] (\hhh*2,-\vvv*2) rectangle (\hhh*3,-\vvv*1) node[pos=.5] {\small $4$};
\end{tikzpicture} 
\hskip 1mm
\begin{tikzpicture}[baseline=0mm]
\def \hhh{4mm}
\def \vvv{5mm}
\draw[-] (\hhh*0,\vvv*0) rectangle (\hhh*1,\vvv*1) node[pos=.5] {\small $1$};
\draw[-] (\hhh*1,\vvv*0) rectangle (\hhh*2,\vvv*1) node[pos=.5] {\small $1$};
\draw[-] (\hhh*2,\vvv*0) rectangle (\hhh*3,\vvv*1) node[pos=.5] {\small $2'$};
\draw[-] (\hhh*3,\vvv*0) rectangle (\hhh*4,\vvv*1) node[pos=.5] {\small $2$};
\draw[-] (\hhh*4,\vvv*0) rectangle (\hhh*5,\vvv*1) node[pos=.5] {\small $3$};
\draw[-] (\hhh*1,-\vvv*1) rectangle (\hhh*2,\vvv*0) node[pos=.5] {\small $2$};
\draw[-] (\hhh*2,-\vvv*1) rectangle (\hhh*3,\vvv*0) node[pos=.5] {\small $3$};
\draw[-] (\hhh*2,-\vvv*2) rectangle (\hhh*3,-\vvv*1) node[pos=.5] {\small $4$};
\end{tikzpicture} 
\hskip 1mm
\begin{tikzpicture}[baseline=0mm]
\def \hhh{4mm}
\def \vvv{5mm}
\draw[-] (\hhh*0,\vvv*0) rectangle (\hhh*1,\vvv*1) node[pos=.5] {\small $1$};
\draw[-] (\hhh*1,\vvv*0) rectangle (\hhh*2,\vvv*1) node[pos=.5] {\small $1$};
\draw[-] (\hhh*2,\vvv*0) rectangle (\hhh*3,\vvv*1) node[pos=.5] {\small $2'$};
\draw[-] (\hhh*3,\vvv*0) rectangle (\hhh*4,\vvv*1) node[pos=.5] {\small $2$};
\draw[-] (\hhh*4,\vvv*0) rectangle (\hhh*5,\vvv*1) node[pos=.5] {\small $3'$};
\draw[-] (\hhh*1,-\vvv*1) rectangle (\hhh*2,\vvv*0) node[pos=.5] {\small $2$};
\draw[-] (\hhh*2,-\vvv*1) rectangle (\hhh*3,\vvv*0) node[pos=.5] {\small $3$};
\draw[-] (\hhh*2,-\vvv*2) rectangle (\hhh*3,-\vvv*1) node[pos=.5] {\small $4$};
\end{tikzpicture} 
\hskip 1mm
\begin{tikzpicture}[baseline=0mm]
\def \hhh{4mm}
\def \vvv{5mm}
\draw[-] (\hhh*0,\vvv*0) rectangle (\hhh*1,\vvv*1) node[pos=.5] {\small $1$};
\draw[-] (\hhh*1,\vvv*0) rectangle (\hhh*2,\vvv*1) node[pos=.5] {\small $1$};
\draw[-] (\hhh*2,\vvv*0) rectangle (\hhh*3,\vvv*1) node[pos=.5] {\small $2$};
\draw[-] (\hhh*3,\vvv*0) rectangle (\hhh*4,\vvv*1) node[pos=.5] {\small $2$};
\draw[-] (\hhh*4,\vvv*0) rectangle (\hhh*5,\vvv*1) node[pos=.5] {\small $3'$};
\draw[-] (\hhh*1,-\vvv*1) rectangle (\hhh*2,\vvv*0) node[pos=.5] {\small $2$};
\draw[-] (\hhh*2,-\vvv*1) rectangle (\hhh*3,\vvv*0) node[pos=.5] {\small $3$};
\draw[-] (\hhh*2,-\vvv*2) rectangle (\hhh*3,-\vvv*1) node[pos=.5] {\small $4$};
\end{tikzpicture} 
\hskip 1mm
\begin{tikzpicture}[baseline=0mm]
\def \hhh{4mm}
\def \vvv{5mm}
\draw[-] (\hhh*0,\vvv*0) rectangle (\hhh*1,\vvv*1) node[pos=.5] {\small $1$};
\draw[-] (\hhh*1,\vvv*0) rectangle (\hhh*2,\vvv*1) node[pos=.5] {\small $1$};
\draw[-] (\hhh*2,\vvv*0) rectangle (\hhh*3,\vvv*1) node[pos=.5] {\small $2$};
\draw[-] (\hhh*3,\vvv*0) rectangle (\hhh*4,\vvv*1) node[pos=.5] {\small $2$};
\draw[-] (\hhh*4,\vvv*0) rectangle (\hhh*5,\vvv*1) node[pos=.5] {\small $3$};
\draw[-] (\hhh*1,-\vvv*1) rectangle (\hhh*2,\vvv*0) node[pos=.5] {\small $2$};
\draw[-] (\hhh*2,-\vvv*1) rectangle (\hhh*3,\vvv*0) node[pos=.5] {\small $3$};
\draw[-] (\hhh*2,-\vvv*2) rectangle (\hhh*3,-\vvv*1) node[pos=.5] {\small $4$};
\end{tikzpicture} 
\begin{tikzpicture}[baseline=0mm]
\def \hhh{4mm}
\def \vvv{5mm}
\node[] at (\hhh*1,-\vvv*0.5) {(when $\nu=(5,2,1)$)};
\end{tikzpicture} 
\vskip 3mm
\noindent
\begin{tikzpicture}[baseline=0mm]
\def \hhh{4mm}
\def \vvv{5mm}
\draw[-] (\hhh*0,\vvv*0) rectangle (\hhh*1,\vvv*1) node[pos=.5] {\small $1$};
\draw[-] (\hhh*1,\vvv*0) rectangle (\hhh*2,\vvv*1) node[pos=.5] {\small $1$};
\draw[-] (\hhh*2,\vvv*0) rectangle (\hhh*3,\vvv*1) node[pos=.5] {\small $2'$};
\draw[-] (\hhh*3,\vvv*0) rectangle (\hhh*4,\vvv*1) node[pos=.5] {\small $3$};
\draw[-] (\hhh*1,-\vvv*1) rectangle (\hhh*2,\vvv*0) node[pos=.5] {\small $2$};
\draw[-] (\hhh*2,-\vvv*1) rectangle (\hhh*3,\vvv*0) node[pos=.5] {\small $2$};
\draw[-] (\hhh*3,-\vvv*1) rectangle (\hhh*4,\vvv*0) node[pos=.5] {\small $4$};
\draw[-] (\hhh*2,-\vvv*2) rectangle (\hhh*3,-\vvv*1) node[pos=.5] {\small $3$};
\end{tikzpicture} 
\hskip 5mm
\begin{tikzpicture}[baseline=0mm]
\def \hhh{4mm}
\def \vvv{5mm}
\draw[-] (\hhh*0,\vvv*0) rectangle (\hhh*1,\vvv*1) node[pos=.5] {\small $1$};
\draw[-] (\hhh*1,\vvv*0) rectangle (\hhh*2,\vvv*1) node[pos=.5] {\small $1$};
\draw[-] (\hhh*2,\vvv*0) rectangle (\hhh*3,\vvv*1) node[pos=.5] {\small $2'$};
\draw[-] (\hhh*3,\vvv*0) rectangle (\hhh*4,\vvv*1) node[pos=.5] {\small $3'$};
\draw[-] (\hhh*1,-\vvv*1) rectangle (\hhh*2,\vvv*0) node[pos=.5] {\small $2$};
\draw[-] (\hhh*2,-\vvv*1) rectangle (\hhh*3,\vvv*0) node[pos=.5] {\small $2$};
\draw[-] (\hhh*3,-\vvv*1) rectangle (\hhh*4,\vvv*0) node[pos=.5] {\small $3$};
\draw[-] (\hhh*2,-\vvv*2) rectangle (\hhh*3,-\vvv*1) node[pos=.5] {\small $4$};
\end{tikzpicture} 
\hskip 5mm
\begin{tikzpicture}[baseline=0mm]
\def \hhh{4mm}
\def \vvv{5mm}
\draw[-] (\hhh*0,\vvv*0) rectangle (\hhh*1,\vvv*1) node[pos=.5] {\small $1$};
\draw[-] (\hhh*1,\vvv*0) rectangle (\hhh*2,\vvv*1) node[pos=.5] {\small $1$};
\draw[-] (\hhh*2,\vvv*0) rectangle (\hhh*3,\vvv*1) node[pos=.5] {\small $2'$};
\draw[-] (\hhh*3,\vvv*0) rectangle (\hhh*4,\vvv*1) node[pos=.5] {\small $2$};
\draw[-] (\hhh*1,-\vvv*1) rectangle (\hhh*2,\vvv*0) node[pos=.5] {\small $2$};
\draw[-] (\hhh*2,-\vvv*1) rectangle (\hhh*3,\vvv*0) node[pos=.5] {\small $3$};
\draw[-] (\hhh*3,-\vvv*1) rectangle (\hhh*4,\vvv*0) node[pos=.5] {\small $3$};
\draw[-] (\hhh*2,-\vvv*2) rectangle (\hhh*3,-\vvv*1) node[pos=.5] {\small $4$};
\end{tikzpicture} 
\hskip 5mm
\begin{tikzpicture}[baseline=0mm]
\def \hhh{4mm}
\def \vvv{5mm}
\draw[-] (\hhh*0,\vvv*0) rectangle (\hhh*1,\vvv*1) node[pos=.5] {\small $1$};
\draw[-] (\hhh*1,\vvv*0) rectangle (\hhh*2,\vvv*1) node[pos=.5] {\small $1$};
\draw[-] (\hhh*2,\vvv*0) rectangle (\hhh*3,\vvv*1) node[pos=.5] {\small $2$};
\draw[-] (\hhh*3,\vvv*0) rectangle (\hhh*4,\vvv*1) node[pos=.5] {\small $2$};
\draw[-] (\hhh*1,-\vvv*1) rectangle (\hhh*2,\vvv*0) node[pos=.5] {\small $2$};
\draw[-] (\hhh*2,-\vvv*1) rectangle (\hhh*3,\vvv*0) node[pos=.5] {\small $3$};
\draw[-] (\hhh*3,-\vvv*1) rectangle (\hhh*4,\vvv*0) node[pos=.5] {\small $3$};
\draw[-] (\hhh*2,-\vvv*2) rectangle (\hhh*3,-\vvv*1) node[pos=.5] {\small $4$};
\end{tikzpicture} 
\hskip 3mm
\begin{tikzpicture}[baseline=0mm]
\def \hhh{4mm}
\def \vvv{5mm}
\node[right] at (\hhh*1,-\vvv*0.5) {(when $\nu=(5,3)$).};
\end{tikzpicture} 
\vskip 4mm \noindent
On the other hand, $P_{\la/\mu} = 2 P_{(6,2)} + 6P_{(5,3)} + 6P_{(5,2,1)} + 8P_{(4,3,1)}$.
\end{exm}
\vskip 20mm
 


\begin{thebibliography}{HK}
\bibitem{AS}
F. Ardila, L. Serrano, 
{\it Staircase skew Schur functions are Schur $P$-positive},
J. Algebraic Combin. {\bf 36} (2012) 409--423.

\bibitem{BSS}
G. Benkart, F. Sottile, J. Stroomer,
{\it Tableau switching: algorithms and applications},
J. Combin. Theory Ser. A {\bf 76} (1996), 11--43.

\bibitem{Cho}
S. Cho, 
{\it Littlewood-Richardson rule for Schur P-functions}, Trans. Amer. Math. Soc. {\bf 365} (2013) 939--972.

\bibitem{CNO14}
S.-I. Choi, S.-Y. Nam, Y.-T. Oh,
{\it Bijections among combinatorial models for shifted Littlewood-Richardson coefficients}, J. Combin. Theory Ser. A {\bf 128} (2014), 56--83.

\bibitem{Dew}
E. Dewitt, 
{\it Identities relating Schur $s$-functions and $Q$-functions},
Ph.D. Thesis, University of Michigan, (2012).

\bibitem{GJKKK14}
D. Grantcharov, J. H. Jung, S.-J. Kang, M. Kashiwara, M. H. Kim,
{\it Crystal bases for the quantum queer superalgebra and semistandard decomposition tableaux},
Trans. Amer. Math. Soc.  {\bf 366} no. 1 (2014), 457--489.

\bibitem{GJKKK15}
D. Grantcharov, J. H. Jung, S.-J. Kang, M. Kashiwara, M. H. Kim,
{\it Crystal bases for the quantum queer superalgebra},
J. Eur. Math. Soc. {\bf 17} (2015), 1593--1627.

\bibitem{Kas95}
M. Kashiwara, 
{\it On crystal bases}, 
Representations of groups,
CMS Conf. Proc., 16, Amer. Math. Soc., Providence, RI, (1995) 155--197.

\bibitem{KashNaka}
M. Kashiwara, T. Nakashima, 
{\it Crystal graphs for representations of the $q$-analogue of classical Lie algebras}, J. Algebra {\bf 165} (1994) 295--345.

\bibitem{Mac}
I. G. Macdonald, 
{\it Symmetric functions and Hall polynomials},
Oxford Mathematical Monographs, second ed., Clarendon
Press, Oxford University Press, New York, 1995.

\bibitem{Sa}
B E. Sagan, 
{\it Shifted tableaux, Schur $Q$-functions and a conjecture of R. Stanley},
J. Combin. Theory Ser.A {\bf 45} (1987), 62--103.

\bibitem{Sch}
I. Schur, {\it \"{U}ber die darstellung der symmetrischen und der alternierenden gruppe durch gebrochene lineare substitutionen}, J. Reine Angew. Math. {\bf 139} (1911), 155--250.

\bibitem{Srgv85}
A. N. Sergeev, {\it Tensor algebra of the identity representation as a module over the Lie superalgebras $Gl(n,m)$ and $Q(n)$}, 
Mat. Sb. (N.S.) {\bf 123(165)} (1984), 422--430 (Russian).

\bibitem{Se}
L. Serrano, 
{\it The shifted plactic monoid}, 
Math. Z. {\bf 266} (2010), 363--392.

\bibitem{SW}
K. M. Shaw, S. van Willigenburg, 
{\it Multiplicity free expansions of Schur $P$-functions}, 
Ann. Combin. {\bf 11} (2007), 69--77.

\bibitem{Ste}
J. R. Stembridge, 
{\it Shifted tableaux and the projective representations of symmetric groups}, Adv. Math. {\bf 74} (1989), 87--134.

\bibitem{Wo}
D. R. Worley,
{\it A theory of shifted Young tableaux}, Ph.D. thesis, MIT, 1984.
\end{thebibliography}
\end{document}